\documentclass[11pt]{elsarticle}

\usepackage{amssymb}
\usepackage{xcolor}
\usepackage{mathrsfs}
\usepackage{amscd}
\usepackage{appendix}
\usepackage{geometry}
\usepackage{setspace}

\usepackage{amsfonts}
\usepackage{amsmath}
\newtheorem{theorem}{Theorem}[section] 
\newtheorem{lemma}[theorem]{Lemma}     

\newtheorem{remark}[theorem]{Remark}

\newtheorem{Definition}[theorem]{Definition}

\numberwithin{equation}{section}

\journal{arXiv}
\begin{document}
\newgeometry{left=2.5cm,right=2.5cm,top=2.2cm,bottom=2.5cm}
\begin{frontmatter}



\title{Invariant Curves of Almost Periodic Reversible Mappings}


\author[1]{Daxiong Piao}\fnref{b}
\fntext[b]{Corresponding author}
\ead{dxpiao@ouc.edu.cn}

\author[1,2]{Xinli Zhang}
\ead{zxl@qust.edu.cn}


\address[1]{School of Mathematical Sciences, Ocean University of China,Qingdao 266100, P. R. China}
\address[2]{School of Mathematics and Physics, Qingdao University of Science and Technology, Qingdao, 266061, P. R. China}

\begin{abstract}
In this paper, we prove some invariant curve theorems for the planar almost periodic reversible mappings. As an application, we will discuss the existence of almost
periodic solutions and the boundedness of all solutions for the nonlinear oscillator
$x''+g(x)x'+\varpi^{2}x+\varphi(x)=f(t)$
with $f(t)$ almost periodic.
\end{abstract}


\end{frontmatter}

\begin{spacing}{1.25}

\section{Introduction}
In this paper, we investigate the existence of invariant curves of the following planar mapping
 \begin{eqnarray}\label{equ11}
\mathfrak{M}:\ \ \left\{\begin{array}{ll}x_{1}=x+\alpha+y+f(x,y),\\
 y_{1}=y+g(x,y),
 \end{array}\right.\ \ \ \ \ \ (x,y)\in \mathbb{R}\times[a,b],
\end{eqnarray}
where $f$ and $g$ are almost periodic in $x$ with the frequency $\omega=(\cdots, \omega_{\lambda}, \cdots)_{\lambda\in\mathbb{Z}}$ and admits a rapidly converging Fourier series expansion. $\alpha$, $a<b $ are positive constants. And we ask under what assumptions the mapping $\mathfrak{M}$ has invariant curves.

If $f$ and $g$ are real analytic, sufficiently small and quasi-periodic  in $x$ with the frequency $\omega=(\omega_1,\omega_2 \cdots, \omega_m)$ which is sufficiently
`incommensurable' together with $2\pi\alpha^{-1}$, and the $\mathfrak{M}$ in (\ref{equ11}) is an exact symplectic map, Zharnitsky \cite{Zh} proves
the existence of invariant curves and applies it to answer a question asked by Levi and Zehnder \cite{LZ}, that is the boundedness of solutions of the
Fermi-Ulam model. His proof is based on the Lagrangian approach introduced by Moser \cite{Mo0} and used by Levi and Moser in \cite{LM} to present a proof of the twist theorem.

Instead of the exact symplecticity of $\mathfrak{M}$ in (\ref{equ11}), Liu \cite{LI} suppose that $\mathfrak{M}$ is
reversible with respect to the involution $\Psi: (x, y) \rightarrow (-x, y)$, that is, $\Psi\mathfrak{M}\Psi = \mathfrak{M}^{-1}$. When $f$ and $g$ are real analytic, sufficiently small and quasi-periodic in $x$ and the frequency $\omega$ satisfies the Diophantine condition,that is $\left|\langle k,\omega\rangle-2\pi j\alpha^{-1}\right|\geq\frac{c_{0}}{|k|^{\sigma_0}},\  (k,j)\in\mathbb{Z}^m\times\mathbb{Z}\backslash \{0,0\}$, for some positive constants $c_0, \sigma_0$, Liu  proved that the mapping $\mathfrak{M}$ has an invariant curve. As an application, he proved the existence of quasi-periodic solutions and the boundedness of
solutions for a pendulum-type equation and an asymmetric oscillator.

How to extend the invariant curve theorems of planar twist mappings with quasi-periodic perturbations to almost periodic perturbations is a longstanding problem.
The difficulty still comes from the so called `small divisor problem'. In the case of quasi-periodicity, the proper Diophantine condition on the frequency can guarantee the KAM iteration convergent. But in the case of almost periodicity, because the dimension of the frequency is infinite, corresponding to Diophantine condition, it seems difficult to find a method to describe the nonresonance condition for the frequency. Recently, Huang, Li and Liu have made a breakthrough in this problem (see \cite{HLL,HLL2}). They, using the `spatial structure' and `weight function' introduced by P\"{o}schel \cite{P} and some elaborate techniques, succussed in proving some invariant curve theorems for planar twist mappings with some kind of almost periodic perturbations.

In \cite{HLL}, Huang-Li-Liu considered a special case of the mapping $\mathfrak{M}$ in (\ref{equ11})
\begin{eqnarray*}
\mathfrak{M_{1}}:\ \ \left\{\begin{array}{ll}x_{1}=x+y+f(x,y),\\
 y_{1}=y+g(x,y),
 \end{array}\right.\ \ \ \ \ \ (x,y)\in \mathbb{R}\times[a,b],
\end{eqnarray*}
where the perturbations $f$ and $g$ are almost periodic in $x$ with the frequency $\omega=(\cdots, \omega_{\lambda}, \cdots)_{\lambda\in\mathbb{Z}}$
When $f$ and $g$ are real analytic and small enough, $\mathfrak{M_{1}}$ possesses the intersection property, the frequency $\omega$ satisfies some kind nonresonance condition, they obtained the invariant curve theorem.

Furthermore, Huang-Li-Liu \cite{HLL2} studied the small twist mapping
\begin{eqnarray}\label{Mdelta}
\mathfrak{M_{\delta}}:\ \ \left\{\begin{array}{ll}x_{1}=x+\beta+\delta y+f(x,y,\delta),\\
 y_{1}=y+g(x,y,\delta),
 \end{array}\right.\ \ \ \ \ \ (x,y)\in \mathbb{R}\times[a,b],
\end{eqnarray}
where the functions $f$ and $g$ are real analytic and almost periodic in $x$ with the frequency $\omega=(\cdots, \omega_{\lambda}, \cdots)_{\lambda\in\mathbb{Z}},$ $\beta$ is a constant, $0<\delta<1$ is a small parameter. When the mapping $\mathfrak{M_{\delta}}$ in (\ref{Mdelta}) has intersection property, the frequency $\omega$ has some nonresonance condition, they proved the existence of invariant curve for the  mapping $\mathfrak{M_{\delta}}$.

Huang-Li-Liu \cite{HLL2} also investigated  the more general small twist mapping

\begin{eqnarray}\label{Mdelta1}
\mathcal{M}_{\delta}:\ \ \left\{\begin{array}{ll}x_{1}=x+\beta+\delta l(x,y)+\delta f(x,y,\delta),\\
 y_{1}=y+\delta m(x,y)+\delta g(x,y,\delta),
 \end{array}\right.\ \ \ \ \ \ (x,y)\in \mathbb{R}\times[a,b],
\end{eqnarray}
where the functions $l, m, f, g$ are real analytic and almost periodic in $x$ with the frequency  $\omega=(\cdots, \omega_{\lambda}, \cdots)_{\lambda\in\mathbb{Z}},$ $f(x, y, 0) =g(x, y, 0) =0$, $0<\delta<1$ . When  $\omega $ and $2\pi/\beta$ are rationally independent, and $\lim_{T\rightarrow \infty}\int_0^T \frac{\partial l}{\partial x}(x,y)dy\neq 0 $, they proved the existence theorem of invariant curve for the mapping $\mathcal{M}_{\delta}$ for sufficiently small $\delta$. Using this theorem, they obtained infinitely many almost periodic solutions and boundedness of all solutions of the asymmetric ocillator
\begin{eqnarray*}
x''+ax^+ - bx^-=f(t)
\end{eqnarray*}
where $f$ is real analytic and almost periodic in $t$ with the frequency  $\omega=(\cdots, \omega_{\lambda}, \cdots)_{\lambda\in\mathbb{Z}}$.

The works of Huang-Li-Liu \cite{HLL,HLL2} widely broadened the applications of the twist mapping theory to the Littewood's boundedness problems for the almost periodic oscillators.

But if we would like to study the boundedness for the equation \footnote[1]{Kunze-Kupper-Liu \cite{KKL} studied the boundedness for Equation (\ref{equ14}) where $ f(t)$ is quasi-periodic forcing. }
\begin{equation}\label{equ14}
x''+g(x)x'+\varpi^{2}x+\varphi(x)=f(t),
\end{equation}
where $\omega>0$, $g,\varphi$ and $f$ are odd functions, $f(t)$ is a real analytic almost periodic function with the frequency $\omega=(\cdots, \omega_{\lambda}, \cdots)_{\lambda\in\mathbb{Z}}$ , we will find that some invariant curve theorems for the almost periodic reversible mappings should be established.

In this paper, motivated by \cite{HLL},\cite{HLL2} and \cite{LI}, we are going to study the existence problems of the invariant curves for the mappings  $\mathfrak{M}$ , $\mathfrak{M_{\delta}}$ and $\mathcal{M_{\delta}}$, when they are reversible and the perturbations are almost periodic. And then, as an application, we shall consider the boundedness of solutions of Equation (\ref{equ14}).

The early contribution on the invariant curve theorem of reversible systems is due to Moser \cite{Mo1}. He studied the existence of invariant tori of a reversible system depending quasi-periodically on time. Then the results was extende by Moser \cite{Mo2} and Sevryuk \cite{Se}.

The rest of the paper is organized as follows. In Section 2, we state the main invariant curve theorem (Theorem \ref{Th1} ) for the almost periodic reversible mapping (\ref{equ11}). The proof of Theorem \ref{Th1} is given in Section 3,4. The small twist theorems are given in Section 5. In Section 6, we will prove the existence of almost periodic solutions for the nonlinear oscillator (\ref{equ14}) with an almost periodic forcing.

\section{Main results}\label{s2}

We will study the existence of invariant curves of the planar almost periodic mapping with a rapidly converging Fourier series expansion. Firstly,  we introduce some definitions and notations .
\begin{Definition}[\cite{SM}] Assume that $\omega=(\omega_{1},\omega_{2},\cdots,\omega_{d})$ are rationally independent. A function $f:\mathbb{R}\rightarrow\mathbb{R}$ is called real analytic quasi-periodic with frequencies $\omega=(\omega_{1},\omega_{2},\cdots,\omega_{d})$ if there exists a real analytic periodic function
$$F:\theta=(\theta_{1},\theta_{2},\cdots,\theta_{d})\in\mathbb{R^{d}}\rightarrow\mathbb{R},$$
in $\theta_{1},\theta_{2},\cdots,\theta_{d}$ with period $2\pi$ and bounded in a complex neighborhood $\Pi_{r}^{d}=\{(\theta_{1}, \theta_{2},\ldots, \theta_{d})\in\mathbb{C}^{d}:|{\rm Im}\theta_{j}\leq r|,j=1,2,\cdots,d\}$ of $\mathbb{R}^{d}$ for some $r>0$, such that
\begin{eqnarray*}
f(t)=F(\omega_{1}t,\ \ \omega_{2}t,\ \ \ldots, \ \ \omega_{d}t),\ \ \ \ \ \ \forall t\in\mathbb{R}.
\end{eqnarray*}
Here, we call $F(\theta)$ the shell function of $f(t)$.
\end{Definition}

Suppose $F(\theta)$ admits a Fourier series expansion
$$F(\theta)=\sum_{k\in\mathbb{Z}^{d}}f_{k}e^{i\langle k,\theta\rangle},$$
where $k=(k_{1}, k_{2},\ldots, k_{d})$, $k_{j}$ range over all integers and the coefficients $f_{k}$ decay exponentially with $|k|=|k_{1}|+|k_{2}|+\cdots+|k_{d}|,$ then $f(t)$ can be represented as a Fourier series of the type from the definition,
$$f(t)=\sum_{k\in\mathbb{Z}^{d}}f_{k}e^{i\langle k,\omega\rangle t}.$$
\begin{Definition}[\cite{SM}] For $r>0$, let $Q_{r}(\omega)$ be the set of real analytic quasi-periodic functions $f$ with frequency $\omega$ such that the shell functions $F$ are bounded on the subset $\Pi_{r}^{d}$ with the supremum norm
$$|F|_{r}=\displaystyle\sup_{\theta\in\Pi_{r}^{d}}\left|F(\theta)\right|=\displaystyle\sup_{\theta\in\Pi_{r}^{d}}\left|\sum_{k}f_{k}e^{i\langle k,\theta\rangle}\right|<+\infty.$$
Thus we define $|f|_{r}:=|F|_{r}.$
\end{Definition}
\begin{Definition}[\cite{Di}] Let $X$ be a complex Banach space. A function $f:U\subseteq X\rightarrow\mathbb{C}$, where $U$
 is an open subset of $X$, is called analytic if $f$ is continuous on $U$, and $f|_{U\cap X_{1}}$ is analytic in the classical sense as
 a function of several complex variables for each finite dimensional subspace $X_{1}$ of $X$.
\end{Definition}


\begin{Definition}[\cite{P}]
Suppose $ \mathcal{S}$ is a family of finite subset of $\mathbb{Z}$. We say that $ \mathcal{S}$ has a \textbf{spatial structure}, If
$$
A,B\in \mathcal{S},\,\, A\cap B \neq \emptyset \Rightarrow A\cup B\in \mathcal{S}.
$$
A nonnegative set function $[ \cdot ] :A\rightarrow [A]$ defined on $ \mathcal{S}\cap\mathcal{S}=\{A\cap B:A,B\in \mathcal{S}\}$ is called an \textbf{weight function} of $\mathcal{S}$ if
\begin{equation*}
A\subseteq B \Rightarrow [A]\leqslant [B],
\end{equation*}
\begin{equation}\label{22}
A\cap B\neq \emptyset \Rightarrow [A\cup B]+[A\cap B]\leqslant[A]+[B].
\end{equation}

For $k\in \mathbb{Z}^\mathbb{Z}$, we define the norm of $k$ as
$$
|k|=\sum_{\lambda\in \mathbb{Z}}|k_{\lambda}|,
$$
and define the support of $k$ as
 $$
 {\rm supp} k=\{\lambda: k_{\lambda}\neq 0\}.
 $$

 We define the weight of the support of $k$ as
 $$
 [[k]]=\min_{{\rm supp} k\subseteqq A \in \mathcal{S}}[A].
 $$

 For $k\in \mathbb{Z}^\mathbb{Z}$ with ${\rm supp} k$ finite and $\theta\in \mathbb{C}^{\mathbb{Z}}$  we define the inner product of $k$ and $\theta$ as
 $$
 \langle k,\theta \rangle=\sum_{\lambda\in \mathbb{Z}} k_{\lambda}\theta_{\lambda}.
 $$
\end{Definition}
\begin{Definition}[\cite{HLL}]\label{2.5} Assume $\omega=(\cdots,\omega_{\lambda},\cdots)_{\lambda\in\mathbb{Z}}$ is a bilateral infinite sequence frequency, its any finite segments are  rationally independent. A function $f:\mathbb{R}\rightarrow\mathbb{R}$ is called real analytic almost periodic with frequency $\omega=(\cdots,\omega_{\lambda},\cdots)\in\mathbb{R^{\mathbb{Z}}}$ if there exists a real analytic periodic function
$$F:\theta=(\cdots,\theta_{\lambda},\cdots)\in\mathbb{R^{\mathbb{Z}}}\rightarrow\mathbb{R},$$
which admit a rapidly converging Fourier series expansion
$$F(\theta)=\sum_{A\in\mathcal{S}}F_{A}(\theta),$$
where
$$F_{A}(\theta)=\sum_{{\rm supp}k\subseteq A}f_{k}e^{i\langle k,\theta\rangle},$$
and
$\mathcal{S}$ has spatial structure with $\mathbb{Z}=\displaystyle{\cup_{A\in\mathcal{S}}A}$, such that
$f(t)=F(\omega t)$ for all $t\in\mathbb{R}$, where $F$ is $2\pi$-periodic in each variable and bounded in a complex neighborhood $\Pi_{r}=\{\theta=(\cdots,\theta_{\lambda},\ldots)\in\mathbb{C}^{\mathbb{Z}}:|{\rm Im}\theta|_{\infty}\leq r\}$ for some $r>0$, where $|{\rm Im}\theta|_{\infty}=\sup_{\lambda\in\mathbb{Z}}|{\rm Im}\theta_{\lambda}|$. Here $F(\theta)$ is called the shell function of $f(t)$.
Thus $f(t)$ can be represented as a Fourier series of the type
\begin{eqnarray}\label{equ21}
f(t)=\sum_{A\in\mathcal{S}}\sum_{{\rm supp}k\subseteq A}f_{k}e^{i\langle k,\omega\rangle t}.
\end{eqnarray}
\end{Definition}

Denote by $AP(\omega)$ the set of all real analytic almost periodic functions with the frequency $\omega$ defined by Definition \ref{2.5}.
\begin{Definition}[\cite{HLL}]\label{lemma26} For $r>0$, let $AP_{r}(\omega)$ be the set of real analytic
almost periodic functions $f$ with the frequency $\omega=(\cdots,\omega_{\lambda},\cdots)$ such that the shell functions $F$ are bounded on the subset $\Pi_{r}$ with the norm
$$||F||_{m,r}:=\sum_{A\in\mathcal{S}}|F_{A}|_{r}e^{m[A]}=\sum_{A\in\mathcal{S}}|f_{A}|_{r}e^{m[A]}<+\infty,$$
where $m>0$ is a constant, $[A]=1+\sum_{i\in A}\log^{\varrho}(1+|i|)$ with $\varrho>2$ is a weight function, and
$$|F_{A}|_{r}=\displaystyle\sup_{\theta\in\Pi_{r}}\left|\sum_{{\rm supp}k\subseteq A}f_{k}e^{i\langle k,\theta\rangle}\right|=|f_{A}|_{r}.$$
Hence we define
$$||f||_{m,r}:=||F||_{m,r}.$$
\end{Definition}

If $f(\cdot, y)\in AP_{r}(\omega)$  and the corresponding shell functions $F(\theta,y)$ are real analytic in the domain
$D(r,s)=\{(\theta, y)\in\mathbb{C}^{\mathbb{Z}}\times\mathbb{C}
:|{\rm Im}\theta|_{\infty}<r, |y|<s\},$
we define
$$||f||_{m,r,s}:=\sum_{A\in\mathcal{S}}|F_{A}|_{r,s}e^{m[A]}=\sum_{A\in\mathcal{S}}|f_{A}|_{r,s}e^{m[A]},$$
where
$$|F_{A}|_{r,s}=\displaystyle\sup_{(\theta,y)\in D{(r,s)}}\left|\sum_{{\rm supp}k\subseteq A}f_{k}e^{i\langle k,\theta\rangle}\right|=|f_{A}|_{r,s}.$$

The real analytic periodic functions have the following properties.
\begin{lemma}[\cite{HLL}]\label{lemma27}
The set $AP(\omega)$ has following properties:

(1) Let $f(t), g(t)\in AP(\omega)$, then $f(t)\pm g(t), g(t+f(t))\in AP(\omega)$;

(2) Let $f(t)\in AP(\omega)$, and $\tau=\beta t+f(t)(\beta+f'>0)$, then the inverse relation is given by $t=\beta^{-1}\tau+g(\tau)$
and $g\in AP(\omega/\beta)$. In particular, if $\beta=1$, then $g\in AP(\omega).$
\end{lemma}

Denote
$$\mathbb{Z}_{0}^{\mathbb{Z}}:=\left\{k=(\cdots,k_{\lambda},\cdots)\in\mathbb{Z}^{\mathbb{Z}}: {\rm supp}k\subseteq A, A\in\mathcal{S}\right\}.$$
Throughout this paper, we assume that the frequency $\omega=(\cdots,\omega_{\lambda},\cdots)$ satisfies the nonresonance conditions
\begin{eqnarray}\label{equ22}
|\langle k,\omega\rangle|\geq\frac{\gamma}{\Delta([[k]]){\Delta(|k|)}},\ \ \ \ \ \ 0\neq k\in\mathbb{Z}_{0}^{\mathbb{Z}},
\end{eqnarray}
where $\gamma$ is a positive constant and $\Delta$ is some fixed approximation function which is defined as follows.

\begin{Definition}[\cite{Ru}] A nondecreasing function $\Delta:[1,\infty)\rightarrow[1,\infty)$ is called an approximation function, if
$$ \frac{\log\Delta(t)}{t}\searrow0,\ \ \ \ 1\leq t\rightarrow\infty,$$
and
$$\int_{1}^{\infty}\frac{\log\Delta(t)}{t^{2}}dt<\infty.$$
\end{Definition}

Now we state our main result.
\begin{theorem}\label{Th1}
Suppose that the almost periodic mapping $\mathfrak{M}$ given by (1.1) is reversible with respect to the involution
$\Psi:(x,y)\mapsto(-x,y)$, that is, $\mathfrak{M}\Psi\mathfrak{M}=\Psi.$ We assume that for every $y$
, $f(\cdot, y), g(\cdot, y)\in AP_{r}(\omega)$, and the corresponding shell
functions $F(\theta, y), G(\theta, y)$ are real analytic in the domain $D(r,s)=\{(\theta, y)\in\mathbb{C}^{\mathbb{Z}}\times\mathbb{C}
:|{\rm Im}\theta|_{\infty}<r, |y|<s\}.$ Furthermore, we assume that
\begin{eqnarray}\label{equ23}
\left|\langle k,\omega\rangle\frac{\alpha}{2\pi}-j\right|\geq\frac{\gamma_{0}}{\Delta([[k]]){\Delta(|k|)}},\ \ \ \ \ k\in\mathbb{Z}_{0}^{\mathbb{Z}}\setminus\{0\}, j\in\mathbb{Z}\setminus\{0\}
\end{eqnarray}
for some positive $\gamma_{0}$. Then there is a positive $\varepsilon_{0}=\varepsilon_{0}(r,s,m,c_{0},\Delta)$
such that if $f,g$ satisfy the following smallness condition
$$\|f\|_{m,r,s}+\|g\|_{m,r,s}<\varepsilon_{0},$$
then the almost periodic mapping $\mathfrak{M}$ has an invariant curve $\Gamma$ and the restriction of $\mathfrak{M}$
onto $\Gamma$ is
$$\mathfrak{M}|_{\Gamma}: x_{1}=x+\alpha.$$
The invariant curve $\Gamma$ is of the form $y=\phi(x)$ with $\phi\in AP_{r'}(\omega)$ for some $r'<r$, and $\|\phi\|_{m',r'}<s, 0<m'<m.$
\end{theorem}　
\begin{remark}
If all the conditions of Theorem \ref{Th1} hold, then the mapping $\mathfrak{M}$ has many invariant curves $\Gamma$ which can be labeled by the form
$$\mathfrak{M}|_{\Gamma}: x_{1}=x+\alpha.$$
of the restriction of $\mathfrak{M}$ onto $\Gamma$.
\end{remark}
\begin{remark}
There is an approximation function $\Delta$ such that for suitable $c_{0}$, the set of $\alpha$ satisfying (\ref{equ23}) has positive measure.
The proof can be found in \cite{HLL}(Theorem 3.4).
\end{remark}
\section{The KAM step}\label{s3}

In this section, we will find a sequence of changes of variables such that the transformed mapping
of $\mathfrak{M}$ will be closer to
 \begin{eqnarray}\label{equ31}
\left\{\begin{array}{ll}x_{1}=x+\alpha+y,\\
 y_{1}=y,
 \end{array}\right.
\end{eqnarray}
than the previous one in the narrower domain. This progress is called the KAM iteration.
\subsection{Construction of the Transformation}
 We will construct a change of variables
\begin{eqnarray}\label{equ32}
{\Phi}:\ \ \ \ \ \ \left\{\begin{array}{ll}x=\xi+\varphi(\xi,\eta),\\
 y=\eta+\psi(\xi,\eta),
 \end{array}\right.
\end{eqnarray}
where $\varphi$ and $\psi$ are real analytic and almost periodic in $\xi$. Under this transformation, the original
mapping $\mathfrak{M}$ is changed into the form
\begin{eqnarray}\label{equ33}
\Phi^{-1}\mathfrak{M}\Phi:\ \ \ \ \ \ \left\{\begin{array}{ll}\xi_{1}=\xi+\alpha+\eta+f_{+}(\xi,\eta),\\
 \eta_{1}=\eta+g_{+}(\xi,\eta),
 \end{array}\right.
\end{eqnarray}
where the functions
\begin{eqnarray}\label{equ34}
\begin{array}{ll}
f_{+}(\xi,\eta)=f(\xi+\varphi,\eta+\psi)+\varphi(\xi,\eta)+\psi(\xi,\eta)-\varphi(\xi+\alpha+\eta+f_{+},\eta+g_{+}),\\
g_{+}(\xi,\eta)=g(\xi+\varphi,\eta+\psi)+\psi(\xi,\eta)-\psi(\xi+\alpha+\eta+f_{+},\eta+g_{+}),
\end{array}
\end{eqnarray}
are real analytic almost periodic functions in $\xi$ defined in a smaller domain $D(r_{+},s_{+})$
and $\|f_{+}\|_{m_{+},r_{+},s_{+}}+\|g_{+}\|_{m_{+},r_{+},s_{+}}$ is smaller than $\|f\|_{m,r,s}+\|g\|_{m,r,s}.$

In the following, $c_{1},c_{2},\cdots$ are positive constants depending on $\gamma_{0},\Delta,\gamma$ only.
We assume that $f(\cdot,y)$ and $g(\cdot,y)$ are real analytic in the domain $D(r,s)$ with $0<r<1, 0<s<\frac{1}{2}$, and for each fixed $y$,
$f(\cdot,y), g(\cdot,y)\in AP_{r}(\omega)$. Moreover, we assume that $\omega=(\cdots, \omega_{\lambda},\cdots)$ and  $\alpha$ satisfy the nonresoance conditions (\ref{equ22}) and (\ref{equ23}).

Let
$$\varepsilon=\|f\|_{m,r,s}+\|g\|_{m,r,s}.$$

We will determine the unknown functions $\varphi, \psi$ from (\ref{equ34}). As one did in Hamiltonian systems, we may solve $\varphi, \psi$ from
the following equations
\begin{eqnarray*}
\begin{array}{ll}
\varphi(\xi+\alpha,\eta)-\varphi(\xi,\eta)=\psi(\xi,\eta)+f(\xi,\eta),\\
\psi(\xi+\alpha,\eta)-\psi(\xi,\eta)=g(\xi,\eta),
\end{array}
\end{eqnarray*}
To ensure that the transformed mapping is reversible with respect to the
same involution $\Psi$, the transformation $\Phi$ should commute with the involution $\Psi$, it follows that
\begin{eqnarray}\label{equ35}
\varphi(-\xi,\eta)=-\varphi(\xi,\eta),\ \ \ \ \ \ \ \ \psi(-\xi,\eta)=\psi(\xi,\eta).
\end{eqnarray}
For this reason, we will find $\varphi,\psi$ from the following modified homological equations
\begin{eqnarray}
\begin{array}{ll}\label{equ36}
\varphi(\xi+\alpha,\eta)-\varphi(\xi,\eta)=F(\xi,\eta),\\
\psi(\xi+\alpha,\eta)-\psi(\xi,\eta)=G(\xi,\eta).
\end{array}
\end{eqnarray}
where
\begin{eqnarray}\label{equ37}
\begin{array}{ll}
F(\xi,\eta)=\displaystyle\frac{1}{2}(\psi(\xi,\eta)+f(\xi,\eta)+\psi(-\xi-\alpha,\eta)+f(-\xi-\alpha,\eta)),\\
G(\xi,\eta)=\displaystyle\frac{1}{2}(g(\xi,\eta)-g(-\xi-\alpha,\eta)).
\end{array}
\end{eqnarray}
It is easy to verify that $F(-\xi-\alpha,\eta)=F(\xi,\eta)$
and $G(-\xi-\alpha,\eta)=-G(\xi,\eta)$.
\subsection{Estimates of the Transformation}
In order to solve $\varphi, \psi$ from (\ref{equ36}), we first solve the following homological equation
\begin{eqnarray}\label{equ38}
l(x+\alpha)-l(x)=h(x),
\end{eqnarray}
where $h\in AP_{r}(\omega)$. One has
\begin{lemma}\label{lemma31}
Suppose that $h\in AP_{r}(\omega)$ and $\omega=(\cdots, \omega_{\lambda}, \cdots)$ satisfying the nonresonance condition (2.2). Then
for any $0<r'<r$, Equation (\ref{equ38}) has a unique solution $l\in AP_{r'}(\omega)$ with $\displaystyle\lim_{T\rightarrow\infty}\frac{1}{T}\int_{0}^{T}l(x)dx=0$ if and only if
\begin{eqnarray}\label{equ39}
\lim_{T\rightarrow\infty}\frac{1}{T}\int_{0}^{T}h(x)dx=0
\end{eqnarray}
In this case, we have the following estimate
\begin{eqnarray}\label{equ310}
\|l\|_{m',r'}\leq\gamma_{0}^{-1}\Lambda(r-r')\Lambda(m-m')\|h\|_{m,r}
\end{eqnarray}
for $0<m'<m, \Lambda(\rho)=\sup_{t\geq0}\Delta(t)e^{-\rho t}.$
Moreover, if $h(-x-\alpha)=h(x)$, then $l$ is odd in $x$; if $h(-x-\alpha)=-h(x)$, then $l$ is even in $x$.
\end{lemma}
{\bf Proof}. The proof of (\ref{equ39}) and (\ref{equ310}) are similar with Lemma 4.1 of \cite{HLL}. Here we only prove the last conclusion.
Because of $h\in AP_{r}(\omega)$ and (\ref{equ39}), $h$ can be represented by
$$h(x)=\sum_{k\neq0,k\in\mathbb{Z}_{0}^{\mathbb{Z}}}h_{k}e^{i\langle k,\omega\rangle x}.$$
Let$$l(x)=\sum_{k\in\mathbb{Z}_{0}^{\mathbb{Z}}}l_{k}e^{i\langle k,\omega\rangle x}.$$
When $h(-x-\alpha)=h(x),$ we have
$$h_{k}=h_{-k}e^{i\langle k,\omega\rangle \alpha},$$
which yields that
\begin{eqnarray*}
l(-x)&=&\sum_{k\neq0,k\in\mathbb{Z}_{0}^{\mathbb{Z}}}\frac{1}{e^{i\langle k,\omega\rangle \alpha}-1}h_{k}e^{-i\langle k,\omega\rangle x}
=\sum_{k\neq0,k\in\mathbb{Z}_{0}^{\mathbb{Z}}}\frac{1}{e^{-i\langle k,\omega\rangle \alpha}-1}h_{-k}e^{i\langle k,\omega\rangle x}\\
&=&-\sum_{k\neq0,k\in\mathbb{Z}_{0}^{\mathbb{Z}}}\frac{1}{e^{i\langle k,\omega\rangle \alpha}-1}h_{k}e^{i\langle k,\omega\rangle x}=-l(x).
\end{eqnarray*}
Similarly, one has that $l(-x)=l(x)$ if $h(-x-\alpha)=-h(x)$, which completes the proof of the lemma.\qed

From Lemma \ref{lemma31}, if we can find the functions $\varphi,\psi$ from (\ref{equ36}), then the transformed mapping $\Phi^{-1}\mathfrak{M}\Phi$
is reversible with respect to the involution  $\Psi:(\xi,\eta)\mapsto(-\xi,\eta)$.

Now we solve the functions $\varphi,\psi$ from (\ref{equ36}) and give the estimates for them. Firstly, we write $F$ and $G$ into the Fourier
series of the type
\begin{eqnarray*}
F=\sum_{k\in\mathbb{Z}_{0}^{\mathbb{Z}}}F_{k}(\eta)e^{i\langle k,\omega\rangle \xi},\ \ \ \ \ \ \ G=\sum_{k\in\mathbb{Z}_{0}^{\mathbb{Z}}}G_{k}(\eta)e^{i\langle k,\omega\rangle \xi}.
\end{eqnarray*}
Now we solve the second equation of (\ref{equ36}). By (\ref{equ37}), one has $G_{0}(\eta)=0$. Hence
\begin{eqnarray}\label{equ311}
\psi(\xi,\eta)=\psi_{0}(\eta)+\sum_{k\neq0, k\in\mathbb{Z}_{0}^{\mathbb{Z}}}\frac{1}{e^{i\langle k,\omega\rangle \alpha}-1}G_{k}(\eta)e^{i\langle k,\omega\rangle \xi}.
\end{eqnarray}
By Lemma \ref{lemma31},
\begin{eqnarray}\label{equ312}
\|\psi(\xi,\eta)-\psi_{0}(\eta)\|_{m-\nu,r-\delta,s}\leq\gamma_{0}^{-1}\Lambda(\delta)\Lambda(\nu)\|G\|_{m,r,s},
\end{eqnarray}
where $0<\nu<m, 0<\delta<r$, $\psi_{0}(\eta)=-f_{0}(\eta)$. From the first equation of (\ref{equ36}), we get
\begin{eqnarray}\label{equ313}
\varphi(\xi,\eta)=\sum_{k\in\mathbb{Z}_{0}^{\mathbb{Z}}}\frac{1}{e^{i\langle k,\omega\rangle \alpha}-1}F_{k}(\eta)e^{i\langle k,\omega\rangle \xi}.
\end{eqnarray}
By Lemma \ref{lemma31},
\begin{eqnarray}\label{equ314}
\|\varphi(\xi,\eta)\|_{m-\nu,r-\delta,s}\leq\gamma_{0}^{-1}\Lambda(\delta)\Lambda(\nu)\|F\|_{m,r,s}
\end{eqnarray}
for $0<\nu<m, 0<\delta<r$. Thus
\begin{eqnarray*}
\|\varphi\|_{m-\nu,r-\delta,s}+\|\psi\|_{m-\nu,r-\delta,s}\leq c_{1}\Lambda(\delta)\Lambda(\nu)(\|F\|_{m,r,s}+\|G\|_{m,r,s}).
\end{eqnarray*}

Similarly, one has
\begin{eqnarray}\label{equ315}
\|\varphi\|_{m-2\nu,r-2\delta,s}+\|\psi\|_{m-2\nu,r-2\delta,s}\leq c_{2}\Lambda^{2}(\delta)\Lambda^{2}(\nu)(\|f\|_{m,r,s}+\|g\|_{m,r,s})=c_{2}\Lambda^{2}(\delta)\Lambda^{2}(\nu)\varepsilon
\end{eqnarray}
and by Cauchy's estimate
\begin{eqnarray}\label{equ316}
\left\|\frac{\partial\varphi}{\partial \xi}\right\|_{m-2\nu,r-3\delta,s}+\left\|\frac{\partial\psi}{\partial \xi}\right\|_{m-2\nu,r-3\delta,s}+
\left\|\frac{\partial\varphi}{\partial \eta}\right\|_{m-2\nu,r-2\delta,s}+\left\|\frac{\partial\psi}{\partial \eta}\right\|_{m-2\nu,r-2\delta,s}\leq c_{2}\Lambda^{2}(\delta)\Lambda^{2}(\nu)\varepsilon(\frac{1}{\rho}+\frac{1}{\delta})\nonumber,\\
\end{eqnarray}
for $0<2\nu<m, 0<3\delta<r$ and $0<\rho<s$.

For $0<m_{+}<m, 0<r_{+}<r<1$ and $0<s_{+}<s<\frac{1}{2}$, let
$$\nu=\frac{1}{10}(m-m_{+}),\ \ \ \ \delta=\frac{1}{10}(r-r_{+}),\ \ \ \ \ \rho=\frac{1}{10}(s-s_{+}).$$
Introduce some domains $D_{j}$ between $D(r_{+}, s_{+})$  and $D(r, s)$ by
$$D_{j}=D(r-j\delta, s-j\rho),\ \ \ \ \ \ \ 0\leq j\leq10.$$
From (\ref{equ314})-(\ref{equ316}), one has
\begin{eqnarray}\label{equ317}
\Phi^{-1}(D_{j+1})\subset D_{j},\ \ \mathfrak{M}(D_{j+1})\subset D_{j},\ \ \ \ \ \ \Phi(D_{j+1})\subset D_{j},
\end{eqnarray}
if
\begin{eqnarray}\label{equ318}
\varepsilon c_{2}\Lambda^{2}(\delta)\Lambda^{2}(\nu)(\frac{1}{\rho}+\frac{1}{\delta})<1.
\end{eqnarray}
So the mapping $\Phi^{-1}\mathfrak{M}\Phi$ is well defined in the domain $D(r_{+},s_{+})$ and maps this domain into $D_{6}$.
\subsection{Estimates of the New Perturbation}
Similar to Subsection 4.3  of \cite{HLL}, one proves that
$f_{+}(\cdot, \eta),g_{+}(\cdot, \eta)\in AP(\omega)$ are well defined by (\ref{equ34}) and are real analytic in $D(r_{+},s_{+})$ if $\varphi(\cdot, \eta),\psi(\cdot, \eta)\in AP(\omega)$. Moreover, one has
$$\|f_{+}\|_{m_{+},r_{+},s_{+}}<6\delta,\ \ \ \ \ \ \|g_{+}\|_{m_{+},r_{+},s_{+}}<6\rho.$$

In the following, we will prove that $\|f_{+}\|_{m_{+},r_{+},s_{+}}+\|g_{+}\|_{m_{+},r_{+},s_{+}}$ is smaller than $\|f\|_{m,r,s}+\|g\|_{m,r,s}.$
By (\ref{equ34}) and (\ref{equ36}), one has
\begin{eqnarray}\label{equ319}
\begin{array}{ll}
f_{+}(\xi,\eta)=f(\xi+\varphi,\eta+\psi)+\varphi(\xi+\alpha,\eta)-\varphi(\xi_{1},\eta_{1})-\frac{1}{2}(G(\xi,\eta)+f(\xi,\eta)+f(-\xi-\alpha,\eta)),\\
g_{+}(\xi,\eta)=g(\xi+\varphi,\eta+\psi)-G(\xi,\eta)+\psi(\xi+\alpha,\eta)-\psi(\xi_{1},\eta_{1}).
\end{array}
\end{eqnarray}
We first estimate the quantity $g(\xi+\varphi,\eta+\psi)-G(\xi,\eta)$. From the reversibility of $\mathfrak{M}$, one has
\begin{eqnarray}\label{equ320}
\begin{array}{ll}
f(-x-\alpha-y-f,y+g)+g(x,y)-f(x,y)=0,\\
g(-x-\alpha-y-f,y+g)+g(x,y)=0.
\end{array}
\end{eqnarray}
Then, for $(\xi,\eta)\in D(r_{+},s_{+})$, one has
\begin{eqnarray}\label{equ321}
|G(\xi,\eta)-g(\xi,\eta)|&=&\left|\frac{1}{2}(g(\xi,\eta)+g(-\xi-\alpha,\eta))\right|\nonumber\\
&=&\left|\frac{1}{2}(g(-\xi-\alpha-\eta-f(\xi,\eta),\eta+g(\xi,\eta))-g(-\xi-\alpha,\eta))\right|\\
&\leq&c_{3}\varepsilon\left(\frac{s_{+}+\varepsilon}{r-r_{+}}+\frac{\varepsilon}{s-s_{+}}\right)\nonumber.
\end{eqnarray}
From (\ref{equ315}), (\ref{equ316}) and (\ref{equ320}), one has
\begin{eqnarray}\label{equ322}
\|g_{+}\|_{m_{+},r_{+},s_{+}}&\leq& c_{4}\varepsilon\Lambda^{2}(\frac{r-r_{+}}{10})\Lambda^{2}(\frac{m-m_{+}}{10})\left(\frac{1}{r-r_{+}}+\frac{1}{s-s_{+}}\right)
(s_{+}+\|f_{+}\|_{m_{+},r_{+},s_{+}}+\|g_{+}\|_{m_{+},r_{+},s_{+}})\nonumber\\&+&c_{3}\varepsilon\left(\frac{s_{+}+\varepsilon}{r-r_{+}}+\frac{\varepsilon}{s-s_{+}}\right)
.
\end{eqnarray}
From the first equation of (\ref{equ320}), one has
\begin{eqnarray*}
\left|\frac{1}{2}(G(\xi,\eta)+f(\xi,\eta)+f(-\xi-\alpha,\eta))-f(\xi,\eta)\right|
\leq c_{3}\varepsilon\left(\frac{s_{+}+\varepsilon}{r-r_{+}}+\frac{\varepsilon}{s-s_{+}}\right).
\end{eqnarray*}
Then, one has
\begin{eqnarray}\label{equ323}
\|f_{+}\|_{m_{+},r_{+},s_{+}}&\leq& c_{4}\varepsilon\Lambda^{2}(\frac{r-r_{+}}{10})\Lambda^{2}(\frac{m-m_{+}}{10})\left(\frac{1}{r-r_{+}}+\frac{1}{s-s_{+}}\right)
(s_{+}+\|f_{+}\|_{m_{+},r_{+},s_{+}}+\|g_{+}\|_{m_{+},r_{+},s_{+}})\nonumber\\&+&c_{3}\varepsilon\left(\frac{s_{+}+\varepsilon}{r-r_{+}}+\frac{\varepsilon}{s-s_{+}}\right)
.
\end{eqnarray}
Now if we choose $\varepsilon$ sufficiently small such that
$$c_{4}\varepsilon\Lambda^{2}(\frac{r-r_{+}}{10})\Lambda^{2}(\frac{m-m_{+}}{10})\left(\frac{1}{r-r_{+}}+\frac{1}{s-s_{+}}\right)<\frac{1}{4},$$
then combining with (\ref{equ322}) and (\ref{equ323}), we have
\begin{eqnarray}\label{equ324}
\|f_{+}\|_{m_{+},r_{+},s_{+}}+\|g_{+}\|_{m_{+},r_{+},s_{+}}\leq c_{5}\varepsilon(s_{+}+\varepsilon)\Lambda^{2}(\frac{r-r_{+}}{10})\Lambda^{2}(\frac{m-m_{+}}{10})
\left(\frac{1}{r-r_{+}}+\frac{1}{s-s_{+}}\right).
\end{eqnarray}
\subsection{The Iteration Lemma}
The above discussions lead to the following lemma.

\begin{lemma}\label{lemma32}

Consider a reversible map
 \begin{eqnarray}\label{equ325}
\mathfrak{M}:\ \ \left\{\begin{array}{ll}x_{1}=x+\alpha+y+f(x,y),\\
 y_{1}=y+g(x,y),
 \end{array}\right.
\end{eqnarray}
 where $f$ and $g$ are real analytic in the domain $D(r,s)$ and almost periodic in $x$ with the frequency $\omega=(\cdots, \omega_{\lambda}, \cdots)_{\lambda\in\mathbb{Z}}$. Assume $\omega, \alpha$ satisfy (\ref{equ22}) and (\ref{equ23}). Let
$$\varepsilon=\|f\|_{m,r,s}+\|g\|_{m,r,s}.$$
Then there is a constant $c_{6}=c_{6}(\gamma_{0},\Delta)>0$ such that for
\begin{eqnarray}\label{equ326}
0<m_{+}<m<1,\ \  0<r_{+}<r<1,\ \ 0<s_{+}<s<1,
\end{eqnarray}
if
\begin{eqnarray}\label{equ327}
\vartheta=c_{6}\varepsilon\Lambda^{2}(\frac{r-r_{+}}{10})\Lambda^{2}(\frac{m-m_{+}}{10})\left(\frac{1}{r-r_{+}}+\frac{1}{s-s_{+}}\right)<\frac{1}{4},\end{eqnarray}
there is a transformation
\begin{eqnarray*}
{\Phi}:\ \ \ \ \ \ \left\{\begin{array}{ll}x=\xi+\varphi(\xi,\eta),\\
 y=\eta+\psi(\xi,\eta),
 \end{array}\right.
\end{eqnarray*}
which is defined in the domain $D(r_{+},s+)$, $\varphi,\psi$ are real analytic and almost periodic in $\xi$.
Under this transformation, the original
mapping $\mathfrak{M}$ is changed into the form
\begin{eqnarray*}
\Phi^{-1}\mathfrak{M}\Phi:\ \ \ \ \ \ \left\{\begin{array}{ll}\xi_{1}=\xi+\alpha+\eta+f_{+}(\xi,\eta),\\
 \eta_{1}=\eta+g_{+}(\xi,\eta),
 \end{array}\right.
\end{eqnarray*}
where the functions $f_{+}$ and $g_{+}$ are
real analytic almost periodic functions in $\xi$ with the frequency $\omega=(\cdots, \omega_{\lambda}, \cdots)_{\lambda\in\mathbb{Z}}$ defined in a smaller domain $D(r_{+},s_{+})$. This map is also reversible with respect to the involution: $(\xi,\eta)\mapsto(-\xi,\eta)$
. Moreover, the following estimates hold:
\begin{eqnarray}\label{equ328}
&&\|\Phi-id\|_{m_{+},r_{+},s_{+}}\leq \vartheta\left(\frac{1}{r-r_{+}}+\frac{1}{s-s_{+}}\right)^{-1},\ \ \ \ \|\partial\Phi-E\|_{m_{+},r_{+},s_{+}}\leq \vartheta,\\
&&\|f_{+}\|_{m_{+},r_{+},s_{+}}+\|g_{+}\|_{m_{+},r_{+},s_{+}}
\leq \vartheta(s_{+}+\varepsilon).
\end{eqnarray}
\end{lemma}
\section{Proof of the Main Results}\label{s4}
In this section, we will use Lemma \ref{lemma32} infinite times to construct a sequence of transformation $\Phi$ to prove Theorem \ref{Th1}.
Denote by $\mathfrak{M}=\mathfrak{M}_{0}$ and restricted to the domain$$D_{0}:|{\rm Im} x|_{\infty}<r_{0},\ \  |y|<s_{0}$$ with $0<r_{0}<1,\ \ 0<s_{0}<1,\ \ 0<m_{0}<1.$
By assumption one has
$$\|f\|_{m_{0},r_{0},s_{0}}+\|g\|_{m_{0},r_{0},s_{0}}<\varepsilon_{0}.$$
Let $$m_{1}=\frac{3}{4}m_{0},\ \ r_{1}=\frac{3}{4}r_{0},\ \ r_{0}=s_{0}=\varepsilon_{0}^{\frac{2}{3}},\ \ s_{1}=\varepsilon_{1}^{\frac{2}{3}},\ \ \varepsilon_{1}=\frac{1}{2}\varepsilon_{0},$$ then we have
$$0<m_{1}<m_{0},\ \  0<r_{1}<r_{0},\ \ 0<s_{1}<s_{0},$$
and
$$\vartheta=c_{6}\varepsilon_{0}\Lambda^{2}(\frac{r_{0}}{40})\Lambda^{2}(\frac{m_{0}}{40})\left(\frac{4}{r_{0}}+\frac{4}{s_{0}}\right)<\frac{1}{4},$$
if $\varepsilon_{0}$ can be chosen sufficiently small. Hence the inequalities (\ref{equ326}),(\ref{equ327}) hold with $n=0$.

Define $$m_{n}=\frac{m_{0}}{2}(1+\frac{1}{2^{n}}),\ \ r_{n}=\frac{r_{0}}{2}(1+\frac{1}{2^{n}}),\ \ s_{n}=\varepsilon_{n}^{\frac{2}{3}},\ \,\varepsilon_{n+1}=\frac{1}{2}\varepsilon_{n},\ \ \  n=0,1,2,\cdots$$ Transforming the mapping $\mathfrak{M}_{0}$
by the coordinate transformation $\Phi=\Phi_{0}$ provided by Lemma \ref{lemma32}, there is a mapping $\mathfrak{M}_{1}=\Phi^{-1}\mathfrak{M}_{0}\Phi$ defined in the domain $D_{1}:|{\rm Im} x|_{\infty}<r_{1},\ \  |y|<s_{1}.$ Applying Lemma 3.2 to the new mapping $\mathfrak{M}_{1}$, there is another coordinate transformation $\Phi_{1}$ and a transformed mapping $\mathfrak{M}_{2}=\Phi_{1}^{-1}\mathfrak{M}_{0}\Phi_{1}$, and proceeding in this way we are led to a sequence of mappings
$$\mathfrak{M}_{n+1}=\Phi_{n}^{-1}\mathfrak{M}_{n}\Phi_{n}\ \ \ \ (n=0,1,2,\cdots)$$
whose domains $D_{n+1}:|{\rm Im} x|_{\infty}<r_{n+1},\ \  |y|<s_{n+1}.$ We proceed by induction and obtaine that the mapping $\mathfrak{M}_{n}$ is well defined in $D_{n}$ and satisfies the appropriate estimate.

By the iteration Lemma \ref{lemma32}, one has a sequence of transformations $\Upsilon_{n}=\Phi_{0}\circ\Phi_{1}\circ\cdots\circ\Phi_{n}$ can be expressed in the form
\begin{eqnarray*}
\left\{\begin{array}{ll}x=\xi+p_{n}(\xi,\eta),\\
 y=\eta+q_{n}(\xi,\eta),
 \end{array}\right.
\end{eqnarray*}
where $$p_{n}=\varphi_{n}+\varphi_{n-1}+\cdots+\varphi_{0},\ \ \ \ q_{n}=\psi_{n}+\psi_{n-1}+\cdots+\psi_{0}.$$
From (\ref{equ328}), it follows that
\begin{eqnarray}\label{equ41}
\|\varphi_{n}\|_{m_{n+1},r_{n+1},s_{n+1}}+\|\psi_{n}\|_{m_{n+1},r_{n+1},s_{n+1}}
\leq(\frac{1}{r_{0}}+\frac{1}{s_{0}})2^{-(n+4)}.
\end{eqnarray}
Then we have
\begin{eqnarray}\label{equ42}
&&\|p_{n}\|_{m_{n+1},r_{n+1},s_{n+1}}+\|q_{n}\|_{m_{n+1},r_{n+1},s_{n+1}}\nonumber\\
&\leq& (\frac{1}{r_{0}}+\frac{1}{s_{0}})\sum_{i=0}^{n}2^{-(i+4)}<(\frac{1}{r_{0}}+\frac{1}{s_{0}})\sum_{i=0}^{\infty}2^{-(i+4)}<\frac{1}{8}(\frac{1}{r_{0}}+\frac{1}{s_{0}}).
\end{eqnarray}
In the following, we will prove
\begin{eqnarray}\label{equ43}
P'_{n1}+P'_{n2}\leq \prod_{i=0}^{n}(1+2\vartheta_{i})-1,
\end{eqnarray}
by induction, where
\begin{eqnarray*}
P'_{n1}&=&\max\left\{\left\|\frac{\partial p_{n}}{\partial \xi}\right\|_{m_{n+1},r_{n+1},s_{n+1}},\ \ \left\|\frac{\partial q_{n}}{\partial \xi}\right\|_{m_{n+1},r _{n+1},s_{n+1}}\right\},\\
P'_{n2}&=&\max\left\{\left\|\frac{\partial p_{n}}{\partial \eta}\right\|_{m_{n+1},r_{n+1},s_{n+1}},\ \ \left\|\frac{\partial q_{n}}{\partial \eta}\right\|_{m_{n+1},r_{n+1},s_{n+1}}\right\},\\
 \vartheta_{i}&=&c_{6}\varepsilon_{i}\Lambda^{2}(\frac{r_{i}-r_{i+1}}{10})\Lambda^{2}(\frac{m_{i}-m_{i+1}}{10})\left(\frac{1}{r_{i}-r_{i+1}}+\frac{1}{s_{i}-s_{i+1}}\right).
\end{eqnarray*}
 For $n=0$, it is obvious by (\ref{equ328}). Let $n\geq1$, from $\Upsilon_{n}=\Upsilon_{n-1}\Phi_{n}$ it follows that
\begin{eqnarray}\label{equ44}
p_{n}=\varphi_{n}+p_{n-1}(\xi+\varphi_{n},\eta+\psi_{n}),\ \ \ \ q_{n}=\psi_{n}+q_{n-1}(\xi+\varphi_{n},\eta+\psi_{n}).
\end{eqnarray}
Then one has
\begin{eqnarray*}
\left\|\frac{\partial p_{n}}{\partial \xi}\right\|_{m_{n+1},r_{n+1},s_{n+1}}<\vartheta_{n}+P'_{(n-1)1}(1+\vartheta_{n})+P'_{(n-1)2}\vartheta_{n},\\
\left\|\frac{\partial q_{n}}{\partial \xi}\right\|_{m_{n+1},r_{n+1},s_{n+1}}<\vartheta_{n}+P'_{(n-1)1}(1+\vartheta_{n})+P'_{(n-1)2}\vartheta_{n},\\
\left\|\frac{\partial p_{n}}{\partial \eta}\right\|_{m_{n+1},r_{n+1},s_{n+1}}<\vartheta_{n}+P'_{(n-1)1}\vartheta_{n}+P'_{(n-1)2}(1+\vartheta_{n}),\\
\left\|\frac{\partial q_{n}}{\partial \eta}\right\|_{m_{n+1},r_{n+1},s_{n+1}}<\vartheta_{n}+P'_{(n-1)1}\vartheta_{n}+P'_{(n-1)2}(1+\vartheta_{n}).
\end{eqnarray*}
Hence
\begin{eqnarray*}
P'_{n1}+P'_{n2}&<&2\vartheta_{n}+(P'_{(n-1)1}+P'_{(n-1)2})(1+2\vartheta_{n})\\
&<&2\vartheta_{n}+(1+2\vartheta_{n})\left\{\prod_{i=0}^{n-1}(1+2\vartheta_{i})-1\right\}<\prod_{i=0}^{n}(1+2\vartheta_{i})-1.
\end{eqnarray*}
Then, by $\vartheta_{i}<\frac{1}{4}$, one has
\begin{eqnarray}\label{equ45}
P'_{n1}+P'_{n2}<\left(\frac{3}{2}\right)^{n}-1
\end{eqnarray}
for all $n$. And (\ref{equ44}) implies
\begin{eqnarray*}
p_{n}-p_{n-1}=\varphi_{n}+p_{n-1}(\xi+\varphi_{n},\eta+\psi_{n})-p_{n-1}(\xi,\eta),\\ q_{n}-q_{n-1}=\psi_{n}+q_{n-1}(\xi+\varphi_{n},\eta+\psi_{n})-q_{n-1}(\xi,\eta).
\end{eqnarray*}
then
\begin{eqnarray*}
\|p_{n}-p_{n-1}\|_{m_{n+1},r_{n+1},s_{n+1}}\leq(\|\varphi_{n}\|_{m_{n+1},r_{n+1},s_{n+1}}+\|\psi_{n}\|_{m_{n+1},r_{n+1},s_{n+1}})(1+P'_{(n-1)1}+P'_{(n-1)2}),\\
\|q_{n}-q_{n-1}\|_{m_{n+1},r_{n+1},s_{n+1}}\leq(\|\varphi_{n}\|_{m_{n+1},r_{n+1},s_{n+1}}+\|\psi_{n}\|_{m_{n+1},r_{n+1},s_{n+1}})(1+P'_{(n-1)1}+P'_{(n-1)2}).
\end{eqnarray*}
By (\ref{equ41}) and (\ref{equ45}), one has
\begin{eqnarray*}
\|p_{n}-p_{n-1}\|_{\frac{m_{0}}{2},\frac{r_{0}}{2}}\leq\left(\frac{1}{r_{0}}+\frac{1}{s_{0}}\right)2^{-4}\left(\frac{3}{4}\right)^{n}<\infty,\\
\|q_{n}-q_{n-1}\|_{\frac{m_{0}}{2},\frac{r_{0}}{2}}\leq\left(\frac{1}{r_{0}}+\frac{1}{s_{0}}\right)2^{-4}\left(\frac{3}{4}\right)^{n}<\infty.
\end{eqnarray*}
Thus $p_{n},q_{n}$ converge to analytic functions of $\xi$ in $|{\rm Im}\xi|<\frac{r_{0}}{2}$ for $\eta=0$. This concludes the proof of the existence of an invariant curve. Hence the proof of Theorem \ref{Th1} is complete.\qed

\section{The small twist theorem}
In this section, we formulate a small twist theorem which is a variant of the invariant curve theorem (Theorem \ref{Th1}) for the almost periodic mapping $\mathfrak{M}$.

In many applications, we may meet the following reversible mappings
 \begin{eqnarray}\label{equ51}
\mathfrak{M}_{\delta}:\ \ \left\{\begin{array}{ll}x_{1}=x+\alpha+\delta y+f(x,y,\delta),\\
 y_{1}=y+g(x,y,\delta),
 \end{array}\right.\ \ \ \ \ \ (x,y)\in \mathbb{R}\times[a,b],
\end{eqnarray}
 where the real analytic functions $f$ and $g$ are almost periodic in $x$ with the frequency $\omega=(\cdots, \omega_{\lambda}, \cdots)_{\lambda\in\mathbb{Z}}$ for each $y$ and admits a rapidly converging Fourier series expansion, $\alpha$ is a positive constant, $0<\delta<1$ is a small parameter.

 We choose the number $\beta$ satisfying the inequalities
\begin{eqnarray}\label{equ52}
\left\{\begin{array}{ll} a+\gamma_{1}\leq\beta\leq b+\gamma_{1},\\
\left|\langle k,\omega\rangle\frac{\alpha+\delta\beta}{2\pi}-j\right|\geq\frac{\gamma_{1}}{\Delta([[k]]){\Delta(|k|)}},\ \ \ \ \ k\in\mathbb{Z}_{0}^{\mathbb{Z}}\setminus\{0\}, j\in\mathbb{Z}\setminus\{0\},
\end{array}\right.
\end{eqnarray}
with some positive constant $\gamma_{1}$, where $\Delta$ is some approximation functions.

\begin{theorem}\label{Th2}
Suppose that the almost periodic mapping $\mathfrak{M}_{\delta}$ given by (\ref{equ51}) is reversible with respect to the involution
$\Psi:(x,y)\mapsto(-x,y)$. We assume that for every $y$
, $f(\cdot, y), g(\cdot, y)\in AP_{r}(\omega)$ with $\omega=(\cdots, \omega_{\lambda}, \cdots)$ satisfying the nonresonance condition (5.2), and the corresponding shell
functions $F(\theta, y), G(\theta, y)$ are real analytic in the domain $D(r,s)=\{(\theta, y)\in\mathbb{C}^{\mathbb{Z}}\times\mathbb{C}
:|{\rm Im}\theta|_{\infty}<r, |y|<s\}.$ There is a positive $\epsilon_{1}>0$
such that if $f,g$ satisfy the following smallness condition
\begin{eqnarray}\label{equ53}
\|f\|_{m,r,s}+\|g\|_{m,r,s}<\delta\epsilon_{1},
\end{eqnarray}
then the almost periodic mapping $\mathfrak{M}_{\delta}$ has an invariant curve $\Gamma$ and the restriction of $\mathfrak{M}_{\delta}$
onto $\Gamma$ is
$$\mathfrak{M}_{\delta}|_{\Gamma}: x_{1}=x+\alpha+\delta\beta.$$
The invariant curve $\Gamma$ is of the form $y=\phi(x)$ with $\phi\in AP_{r'}(\omega)$ for some $r'<r$, and $\|\phi\|_{m',r'}<s, 0<m'<m.$
\end{theorem}　
This is the so called small twist theorem. One can use the same procedure in the proof of Theorem \ref{Th1}
to prove it. We omit it here.
\begin{remark}
If all the conditions of Theorem 5.1 hold, given any $\beta$ satisfying the inequalities (\ref{equ52}), there exists an almost periodic invariant curve $\Gamma$ of $\mathfrak{M}_{\delta}$ with the frequency $\omega=(\cdots, \omega_{\lambda}, \cdots)$ and the restriction of $\mathfrak{M}_{\delta}$ onto
$\Gamma$ is
$$\mathfrak{M}_{\delta}|_{\Gamma}: x_{1}=x+\alpha+\delta\beta.$$
\end{remark}
\begin{remark}
The above conclusion is also true for the following mapping
\begin{eqnarray}\label{equ54}
 \left\{\begin{array}{ll}x_{1}=x+\alpha+\delta h(y)+f(x,y,\delta),\\
 y_{1}=y+g(x,y,\delta),
 \end{array}\right.\ \ \ \ \ \ (x,y)\in \mathbb{R}\times[a,b],
\end{eqnarray}
with $h'(y)\neq0$, if we change the condition (\ref{equ53}) into
\begin{eqnarray*}M(\|f\|_{m,r,s}+\|g\|_{m,r,s})<\epsilon_{1},
\end{eqnarray*}
where $M=\max\{|h'|_{s},1\}.$
If all the conditions of Theorem 5.1 hold, there exists an almost periodic invariant curve $\Gamma$ of $\mathfrak{M}_{\delta}$ with the frequency $\omega=(\cdots, \omega_{\lambda}, \cdots)$ and the restriction of $\mathfrak{M}_{\delta}$ onto
$\Gamma$ is
$$\mathfrak{M}_{\delta}|_{\Gamma}: x_{1}=x+\alpha+\delta\beta.$$
\end{remark}

In the following, we will investigate the mapping
 \begin{eqnarray}\label{equ55}
\mathfrak{M}_{\delta}:\ \ \left\{\begin{array}{ll}x_{1}=x+\alpha+\delta L(x,y)+\delta f(x,y,\delta),\\
 y_{1}=y+\delta M(x,y)+\delta g(x,y,\delta),
 \end{array}\right.\ \ \ \ \ \ (x,y)\in \mathbb{R}\times[a,b],
\end{eqnarray}
 where the real analytic functions $f(\cdot, y), g(\cdot, y)\in AP_{r}(\omega)$ with $\omega=(\cdots, \omega_{\lambda}, \cdots)$ satisfying the nonresonance condition (\ref{equ52}) for every $y$, and admits a rapidly converging Fourier series expansion, $f(x,y,0)=g(x,y,0)=0$, $\alpha$ is a positive constant, $0<\delta<1$ is a small parameter. Suppose that the almost periodic mapping $\mathfrak{M}_{\delta}$ given by (\ref{equ55}) is reversible with respect to the involution
$\Psi:(x,y)\mapsto(-x,y)$. We can obtain the nonresonant and resonant small twist theorems.

\begin{theorem}\label{Th3} In the previous settings, we assume that $\omega_{A}=\{\omega_{\lambda}:\lambda\in A\}$ and $\frac{2\pi}{\alpha}$ are rationally independent for every $A\in\mathcal{S},$ and
$$\lim_{T\rightarrow\infty}\int_{0}^{T}\frac{\partial L}{\partial y}(x,y)dx\neq0.$$
Then there exists $\delta_{0}>0$ such that the mapping $\mathfrak{M}_{\delta}$ has an invariant curve in the domain $\mathbb{R}\times[a,b]$ if $0<\delta<\delta_{0}.$
The invariant curve is of the form $y=\phi(x)$ with $\phi\in AP_{r'}(\omega)$ for some $r'<r$, and $\|\phi\|_{m',r'}<s, 0<m'<m.$
\end{theorem}　
{\bf Proof}. The main idea of the proof is similar to Theorem 3 of Liu \cite{LI} and Theorem 1 of Liu-Song \cite{LS},  so we just give a sketch here.

Since the mapping $\mathfrak{M}_{\delta}$ are reversible with respect to the involution
$\Psi:(x,y)\mapsto(-x,y)$ for all $\delta\in(0,1)$, one has that
\begin{eqnarray*}
L(x,y)+f(x,y,\delta)=L\circ\Psi\circ\mathfrak{M}_{\delta}(x,y)+f(\Psi\circ\mathfrak{M}_{\delta}(x,y),\delta),\\
M(x,y)+g(x,y,\delta)=-M\circ\Psi\circ\mathfrak{M}_{\delta}(x,y)-g(\Psi\circ\mathfrak{M}_{\delta}(x,y),\delta).
\end{eqnarray*}
Let $\delta\rightarrow 0^{+}$, then one has
\begin{eqnarray}\label{equ56}
L(x,y)=L(-x-\alpha,y),\ \ \ \ \ M(x,y)=-M(-x-\alpha,y),\ \ \ \ \ \ \ \ (x,y)\in \mathbb{R}\times[a,b].
\end{eqnarray}
Expanding $L,M$ into Fourier series
\begin{eqnarray*}
L(x,y)=\sum_{A\in\mathcal{S}}\sum_{{\rm supp}k\subseteq A}L_{k}(y)e^{i\langle k,\omega\rangle x},\ \ \ \ \ \  M(x,y)=\sum_{A\in\mathcal{S}}\sum_{{\rm supp}k\subseteq A}M_{k}(y)e^{i\langle k,\omega\rangle x}
\end{eqnarray*}
with $L_{-k}=\overline{L}$ and $M_{-k}=\overline{M}$, then
$$L_{k}=L_{-k}(y)e^{i\langle k,\omega\rangle \alpha},\ \ \ \ \ \ \ M_{k}=-M_{-k}(y)e^{i\langle k,\omega\rangle \alpha}$$
for every $k\in\mathbb{Z}_{0}^{\mathbb{Z}}.$ In particular, we have
$$M_{0}(y)=\lim_{T\rightarrow\infty}\frac{1}{T}\int_{0}^{T}M(x,y)dx\equiv0.$$
From the assumption of Theorem \ref{Th3}, one has
$$L_{0}'(y)=\lim_{T\rightarrow\infty}\int_{0}^{T}\frac{\partial L}{\partial y}(x,y)dx\neq0.$$
Let
\begin{eqnarray*}
H_{1}(x,y)=\sum_{A\in\mathcal{S}}\sum_{\substack{ {\rm supp}k\subseteq A\\
0<\mu [[k]]+\nu |k|<N
}}L_{k}(y)e^{i\langle k,\omega\rangle x},\ \ \ \ \ \  H_{2}(x,y)=\sum_{A\in\mathcal{S}}\sum_{\substack{ {\rm supp}k\subseteq A\\
0<\mu [[k]]+\nu |k|<N
}}M_{k}(y)e^{i\langle k,\omega\rangle x},
\end{eqnarray*}
where $\mu,\nu$ are two small and $N$ is a large positive parameters.

For any $\epsilon_{2}>0$, there exists a positive integer $N=N(L,M,\mu,\nu)$ such that
\begin{eqnarray*}
&&\|L(x,y)-H_{1}(x,y)-L_{0}(y)\|_{m-\mu,r-\nu,s}+\|M(x,y)-H_{2}(x,y)\|_{m-\mu,r-\nu,s}\\
&\leq& e^{-N}(\|L\|_{m,r,s}+\|M\|_{m,r,s})<\epsilon_{2}
\end{eqnarray*}
for $0<\mu<m, 0<\nu<r.$ Moreover, one has
\begin{eqnarray}\label{equ57}
\|L(x,y)-H_{1}(x,y)-L_{0}(y)\|_{m,r,s}+\|M(x,y)-H_{2}(x,y)\|_{m,r,s}<\epsilon_{2}.
\end{eqnarray}
Consider the difference equations
\begin{eqnarray}\label{equ58}
\begin{array}{ll}U(x+\alpha,y)-U(x,y)+H_{1}(x,y)=0,\\
V(x+\alpha,y)-V(x,y)+H_{2}(x,y)=0
\end{array}
\end{eqnarray}
for the unknown functions $U,V$.

From the assumption that $\omega_{A}=\{\omega_{\lambda}:\lambda\in A\}$ and $\frac{2\pi}{\alpha}$ are rationally independent for $A\in\mathcal{S},$ we have
\begin{eqnarray*}
U(x,y)=-\sum_{A\in\mathcal{S}}\sum_{\substack{ {\rm supp}k\subseteq A\\
0<\mu [[k]]+\nu |k|<N
}}\frac{L_{k}(y)}{e^{i\langle k,\omega\rangle \alpha}-1}e^{i\langle k,\omega\rangle x}.\\
V(x,y)=-\sum_{A\in\mathcal{S}}\sum_{\substack{ {\rm supp}k\subseteq A\\
0<\mu [[k]]+\nu |k|<N
}}\frac{M_{k}(y)}{e^{i\langle k,\omega\rangle \alpha}-1}e^{i\langle k,\omega\rangle x}.\\
\end{eqnarray*}
is a solution of (\ref{equ58}). Moreover,
\begin{eqnarray*}
U(-x,y)=-U(x,y),\ \ \ \ V(-x,y)=V(x,y)
\end{eqnarray*}
and there is a positive constant $\varsigma(N)$ such that
\begin{eqnarray*}
\|U\|_{m,r,s}+\|V\|_{m,r,s}\leq\varsigma(\|L\|_{m,r,s}+\|M\|_{m,r,s}).
\end{eqnarray*}
Let
\begin{eqnarray*}
R_{1}(x,y)&=&U(x+\alpha,y)-U(x,y)+L(x,y)-L_{0}(y),\\
R_{2}(x,y)&=&V(x+\alpha,y)-V(x,y)+M(x,y).
\end{eqnarray*}
Then by (\ref{equ57}), one has
\begin{eqnarray}\label{equ59}
\|R_{1}\|_{m,r,s}+\|R_{2}\|_{m,r,s}<\epsilon_{2}.
\end{eqnarray}
Define the change of variables $\mathcal{U}$ by
$$\theta=x+\delta U(x,y),\ \ \ \ \ \ \rho=y+\delta V(x,y).$$
Then the transformed mapping $\mathcal{U}\circ\mathfrak{M}_{\delta}\circ\mathcal{U}^{-1}$ is of the form
\begin{eqnarray}\label{equ510}
\left\{\begin{array}{ll}\theta_{1}=\theta+\alpha+\delta L_{0}(\rho)+\delta\phi_{1}\circ\mathcal{U}^{-1}(\theta,\rho,\delta),\\
\rho_{1}=\rho+\delta\phi_{2}\circ\mathcal{U}^{-1}(\theta,\rho,\delta),
\end{array}\right.
\end{eqnarray}
where
\begin{eqnarray}\label{equ511}
\begin{array}{ll}\phi_{1}(x,y,\delta)=f(x,y,\delta)+R_{1}(x,y)+U(x_{1},y_{1})-U(x+\alpha,y)+L_{0}(y)-L_{0}(y+\delta V(x,y)),\\
\phi_{2}(x,y,\delta)=g(x,y,\delta)+R_{2}(x,y)+V(x_{1},y_{1})-V(x+\alpha,y).
\end{array}
\end{eqnarray}
The functions $\phi_{1},\phi_{2}$ are real analytic and almost periodic with the frequency $\omega=(\cdots, \omega_{\lambda}, \cdots)$ by Lemma \ref{lemma27}.
Hence $\phi_{1}\circ\mathcal{U}^{-1}$ and $\phi_{2}\circ\mathcal{U}^{-1}$ are also real analytic and almost periodic with the frequency $\omega=(\cdots, \omega_{\lambda}, \cdots)$ by Lemma \ref{lemma27}. Similar to \cite{O}, there exists a constant ${M_{1}}>1$ such that
\begin{eqnarray*}
&&\|\phi_{1}\circ\mathcal{U}^{-1}(\cdot,\cdot,\delta)\|_{m,r,s}+\|\phi_{2}\circ\mathcal{U}^{-1}(\cdot,\cdot,\delta)\|_{m,r,s}\\
&\leq&M_{1}\{\|\phi_{1}(\cdot,\cdot,\delta)\|_{m,r,s}+\|\phi_{2}(\cdot,\cdot,\delta)\|_{m,r,s}\}.
\end{eqnarray*}
Since $f(\cdot,\cdot,0)=g(\cdot,\cdot,0)=0$, then there exists a $\delta_{1}>0$ such that
\begin{eqnarray}\label{equ512}
\|f(\cdot,\cdot,\delta)\|_{m,r,s}+\|g(\cdot,\cdot,\delta)\|_{m,r,s}<\epsilon_{3}, \ \ \ \ \ \ \ 0<\delta<\delta_{1},
\end{eqnarray}
where $ 0<\epsilon_{3}=\frac{\epsilon_{1}}{4N_{1}M_{1}}<1, N_{1}=\max\{1,|L'_{0}|_{s}\}.$
Similar to \cite{O}, there is $M_{2}>0$ and $\delta_{2}>0$ , for any $\delta\in(0,\delta_{2})$ such that
\begin{eqnarray}\label{equ513}
&&\|L_{0}(\rho)-L_{0}(y)\|_{m,r,s}<\epsilon_{3}, \\
&&\|U(x_{1},y_{1})-U(x+\alpha,y)\|_{m,r,s}+\|V(x_{1},y_{1})-V(x+\alpha,y)\|_{m,r,s}\leq M_{2}\delta(1+N_{2}),
\end{eqnarray}
where $N_{2}=\|L\|_{m,r,s}+\|M\|_{m,r,s}.$

Choose $$\epsilon_{2}=\frac{\epsilon_{1}}{4N_{1}M_{1}}, \ \ \ \ \ \ \delta_{0}=\min\left\{\delta_{1},\delta_{2},\frac{\epsilon_{1}}{4N_{1}M_{1}M_{2}(1+N_{2})}\right\}.$$
From (\ref{equ59}),(\ref{equ512}) and (\ref{equ513}), one has
$$\|\phi_{1}(\cdot,\cdot,\delta)\|_{m,r,s}+\|\phi_{2}(\cdot,\cdot,\delta)\|_{m,r,s}<\frac{\epsilon_{1}}{4N_{1}M_{1}},$$
which imply that
\begin{eqnarray}\label{equ515}
N_{1}\|\phi_{1}\circ\mathcal{U}^{-1}(\cdot,\cdot,\delta)\|_{m,r,s}+\|\phi_{2}\circ\mathcal{U}^{-1}(\cdot,\cdot,\delta)\|_{m,r,s}<\epsilon_{1}
\end{eqnarray}
for $\delta\in(0,\delta_{0}).$
Hence this mapping $\mathcal{U}\circ\mathfrak{M}_{\delta}\circ\mathcal{U}^{-1}$ meets all assumptions of Remark 5.3 and  has invariant curves.
So undoing the change of variables we obtain the existence of invariant curves of $\mathfrak{M}_{\delta}.$  This end the proof of Theorem \ref{Th3}.\qed

Now we will discuss the resonant case, that is, there is some set $A\in\mathcal{S}$ such that $\omega_{A}=\{\omega_{\lambda}:\lambda\in A\}$ and $\frac{2\pi}{\alpha}$ are rationally dependent. Denote by $\mathbb{S}$ the set of all $A\in\mathcal{S}$ such that there is the integer vector
$k\in\mathbb{Z}^{\mathbb{Z}}_{\mathcal{S}}\setminus\{0\}$ such that ${\rm supp}k\subseteq A, \langle k,\omega\rangle\alpha\in2\pi\mathbb{Z}.$

Then the functions $L$ and $M$ in (\ref{equ55}) can be represented in the form
\begin{eqnarray*}
L(x,y):=\widetilde{L}(x,y)+\widehat{L}(x,y)=\sum_{A\in\mathcal{S}\setminus\mathbb{S}}\sum_{{\rm supp}k\subseteq A}L_{k}(y)e^{i\langle k,\omega\rangle x}+\sum_{A\in\mathbb{S}}\sum_{{\rm supp}k\subseteq A}L_{k}(y)e^{i\langle k,\omega\rangle x},\\
M(x,y):=\widetilde{M}(x,y)+\widehat{M}(x,y)=\sum_{A\in\mathcal{S}\setminus\mathbb{S}}\sum_{{\rm supp}k\subseteq A}M_{k}(y)e^{i\langle k,\omega\rangle x}+\sum_{A\in\mathbb{S}}\sum_{{\rm supp}k\subseteq A}M_{k}(y)e^{i\langle k,\omega\rangle x}.
\end{eqnarray*}
Note that
$$e^{i\langle k,\omega\rangle x}-1\neq0,\ \ \ \ \ {\rm supp}k\subseteq A,\ \ A\in\mathcal{S}\setminus\mathbb{S},$$
and
$$\widehat{L}(x+\alpha,y)\equiv\widehat{L}(x,y),\ \ \ \ \ \widehat{M}(x+\alpha,y)\equiv\widehat{M}(x,y).$$
\begin{theorem}\label{Th4}
In the previous settings, we assume the function $\widehat{L}$ satisfies
$$\widehat{L}(x,y)>0,\ \ \ \ \ \ \frac{\partial \widehat{L}}{\partial y}>0$$
and there is a real analytic function $I(x+\alpha,y)\equiv I(x,y)$ satisfing
\begin{eqnarray}\label{equ516}
&&\frac{\partial I}{\partial y}>0,\\
&&\widehat{L}(x,y)\frac{\partial I}{\partial x}(x,y)+\widehat{M}(x,y)\frac{\partial I}{\partial y}(x,y)\equiv0,
\end{eqnarray}
and two numbers $\widetilde{a},\widetilde{b}$ such that
\begin{eqnarray}\label{equ518}
a<\widetilde{a}&<&\widetilde{b}<b,\nonumber\\
I_{max}(a)<I_{min}(\widetilde{a})\leq I_{max}(\widetilde{a})&<&I_{min}(\widetilde{b})\leq I_{max}(\widetilde{b})<I_{min}(b),
\end{eqnarray}
where
$$I_{min}(y):=min_{x\in\mathbb{R}}I(x,y),\ \ \ \ \ \ I_{max}(y):=max_{x\in\mathbb{R}}I(x,y).$$
Then there exist $\epsilon>0$ and $\widetilde{\delta}_{0}>0$ such that if $0<\delta<\widetilde{\delta}_{0}$ and
\begin{eqnarray}\label{equ519}
\|f(\cdot,\cdot,\delta)\|_{m,r,s}+\|g(\cdot,\cdot,\delta)\|_{m,r,s}<\epsilon,
\end{eqnarray}
and the mapping $\mathfrak{M}_{\delta}$ has an invariant curve in the domain $\mathbb{R}\times[a,b]$ .
The invariant curve is of the form $y=\phi(x)$ with $\phi\in AP_{r'}(\omega)$ for some $r'<r$, and $\|\phi\|_{m',r'}<s, 0<m'<m.$ The constants $\epsilon$ and $\widetilde{\delta}_{0}>0$ depends on $a,b,\widetilde{a},\widetilde{b},L,M,I.$
\end{theorem}　
{\bf Proof}. From (\ref{equ516}) and the definition of the set $\mathbb{S}$, one has
\begin{eqnarray}\label{equ520}
&&\widetilde{L}(-x-\alpha,y)\equiv\widetilde{L}(x,y),\ \ \ \ \ \widetilde{M}(-x-\alpha,y)\equiv-\widetilde{M}(x,y),\\
&&\widehat{L}(-x,y)\equiv\widehat{L}(x,y),\ \ \ \ \ \ \ \ \ \ \widehat{M}(-x,y)\equiv-\widehat{M}(x,y)
\end{eqnarray}
Consider the difference equations
\begin{eqnarray}\label{equ522}
\begin{array}{ll}\widetilde{U}(x+\alpha,y)-\widetilde{U}(x,y)+\widetilde{L}^{N}(x,y)=0,\\
\widetilde{V}(x+\alpha,y)-\widetilde{V}(x,y)+\widetilde{M}^{N}(x,y)=0,
\end{array}
\end{eqnarray}
where
\begin{eqnarray*}
\widetilde{L}^{N}(x,y)=\sum_{A\in\mathcal{S}\setminus\mathbb{S}}\sum_{\substack{ {\rm supp}k\subseteq A\\
0<\mu [[k]]+\nu |k|<N
}}L_{k}(y)e^{i\langle k,\omega\rangle x},\\
 \widetilde{M}^{N}(x,y)=\sum_{A\in\mathcal{S}\setminus\mathbb{S}}\sum_{\substack{ {\rm supp}k\subseteq A\\
0<\mu [[k]]+\nu |k|<N
}}M_{k}(y)e^{i\langle k,\omega\rangle x},
\end{eqnarray*}
with $\mu,\nu$ are two small and $N$ is a large positive parameters.

Since $\langle k,\omega\rangle\alpha\notin2\pi\mathbb{Z}$ for ${\rm supp}k\subseteq A, A\in\mathcal{S}\setminus\mathbb{S},$ one has from the proof of Theorem 5.4 that
\begin{eqnarray*}
\widetilde{U}(x,y)=-\sum_{A\in\mathcal{S}\setminus\mathbb{S}}\sum_{\substack{ {\rm supp}k\subseteq A\\
0<\mu [[k]]+\nu |k|<N
}}\frac{L_{k}(y)}{e^{i\langle k,\omega\rangle \alpha}-1}e^{i\langle k,\omega\rangle x}.\\
\widetilde{V}(x,y)=-\sum_{A\in\mathcal{S}\setminus\mathbb{S}}\sum_{\substack{ {\rm supp}k\subseteq A\\
0<\mu [[k]]+\nu |k|<N
}}\frac{M_{k}(y)}{e^{i\langle k,\omega\rangle \alpha}-1}e^{i\langle k,\omega\rangle x}.\\
\end{eqnarray*}
is a solution of (\ref{equ522}). As we did in the proof of Theorem 5.4, under the transformation  ${\mathcal{U}_{1}}$ by
$$\theta=x+\delta \widetilde{U}(x,y),\ \ \ \ \ \ \rho=y+\delta \widetilde{V}(x,y).$$
Then the transformed mapping ${\mathcal{U}_{1}}\circ\mathfrak{M}_{\delta}\circ{\mathcal{{U}}_{1}}^{-1}$ is of the form
\begin{eqnarray}\label{equ523}
\left\{\begin{array}{ll}\theta_{1}=\theta+\alpha+\delta \widehat{L}(\theta,\rho)+\delta\widetilde{f}(\theta,\rho,\delta),\\
\rho_{1}=\rho+\delta\widehat{M}(\theta,\rho)+\delta\widetilde{g}(\theta,\rho,\delta),
\end{array}\right.
\end{eqnarray}
where the functions $\widetilde{f}$ and $\widetilde{g}$ are very small if $\delta$ is sufficiently small and $N$ is very large. Moreover,
$\widetilde{f}$ and $\widetilde{g}$ are almost periodic in $\theta$ with the frequency $\omega=(\cdots,\omega_{\lambda},\cdots),$  and this mapping
is reversible with respect to the involution $\Psi:(x,y)\mapsto(-x,y).$

In the following, we will construct another transformation ${\mathcal{U}_{2}}$ such that the transformed mapping ${\mathcal{U}_{2}}\circ{\mathcal{U}_{1}}\circ\mathfrak{M}_{\delta}\circ{\mathcal{{U}}_{1}}^{-1}\circ{\mathcal{{U}}_{2}}^{-1}$
is of the form
\begin{eqnarray}\label{equ524}
\left\{\begin{array}{ll}\tau_{1}=\tau+\alpha+\delta \Gamma(\varrho)+\delta\Omega_{1}(\tau,\varrho,\delta),\\
\varrho_{1}=\varrho+\delta\Omega_{2}(\tau,\varrho,\delta).
\end{array}\right.
\end{eqnarray}

For each $\theta\in\mathbb{R}$ and $h\in\mathbb{R}$ with $I(\theta,a)\leq h\leq I(\theta,b),$
define the positive periodic function
\begin{eqnarray*}
\Pi:\ \ \ \ [\overline{I}(a),\underline{I}(b)]\rightarrow\mathbb{R},\ \ \ \ \ \Pi(h)=\int_{0}^{\alpha}\frac{d\theta}{\widehat{L}(\theta,R(\theta,h))},
\end{eqnarray*}
where $R(\theta,h)$ satisfy $I(\theta,R)=h$ and is $\alpha-$periodic in $\theta.$ Moreover,
\begin{eqnarray}\label{equ525}
\Pi'(h)=-\int_{0}^{\alpha}\frac{1}{\widehat{L}^{2}(\theta,R(\theta,h))}\frac{\partial\widehat{L} }{\partial \rho}(\theta,R)\frac{\partial R}{\partial h}d\theta<0.
\end{eqnarray}
Define
\begin{eqnarray*}
K:\ \ \ \ \mathbb{A}\rightarrow\mathbb{R},\ \ \ \ \ \ K(\theta,\rho)=\int_{0}^{\theta}\frac{ds}{\widehat{L}(s,R(s,I(\theta,\rho)))},
\end{eqnarray*}
where $\mathbb{A}=\{(\theta,\rho):\theta\in\mathbb{R}, \widetilde{a}\leq\rho\leq\widetilde{b}\}.$ Moreover, one has
$$K(\theta+\alpha,\rho)=K(\theta,\rho)+\Pi(I(\theta,\rho)),\ \ \ \ \ (\theta,\rho)\in\mathbb{A},$$
and
\begin{eqnarray}\label{equ526}
\widehat{L}(\theta,\rho)\frac{\partial K}{\partial \theta}+\widehat{M}(\theta,\rho)\frac{\partial K}{\partial \rho}=1.
\end{eqnarray}

We define the mapping ${\mathcal{U}_{2}}$ by
\begin{eqnarray}\label{equ527}
\mathcal{U}_{2}:\ \ \ \ \varrho={I}(\theta,\rho),\ \ \ \ \tau=\Gamma(I(\theta,\rho))K(\theta,\rho),
\end{eqnarray}
where $\Gamma(h)=\frac{\alpha}{\Pi(h)}, h\in[\overline{I}(a),\underline{I}(b)].$

We may assume that the function $I(\theta,\rho)$ is even in $\theta$. Otherwise, we use $\frac{{I}(\theta,\rho)+{I}(-\theta,\rho)}{2}$ instead of ${I}(\theta,\rho).$ From (5.17) and (\ref{equ526}), one has
\begin{eqnarray}\label{equ528}
&&\widehat{L}(\theta,\rho)\frac{\partial \tau}{\partial \theta}+\widehat{M}(\theta,\rho)\frac{\partial \tau}{\partial \rho}=\Gamma\circ I,\\
&&\tau(\theta+\alpha,\rho)=\tau(\theta,\rho)+\alpha,\\
&&\frac{\partial \tau}{\partial \theta}(\theta+\alpha,\rho)=\frac{\partial \tau}{\partial \theta}(\theta,\rho),\ \ \ \frac{\partial \tau}{\partial \rho}(\theta+\alpha,\rho)=\frac{\partial \tau}{\partial \rho}(\theta,\rho).
\end{eqnarray}

Therefore the transformed mapping ${\mathcal{U}_{2}}\circ{\mathcal{U}_{1}}\circ\mathfrak{M}_{\delta}\circ{\mathcal{{U}}_{1}}^{-1}\circ{\mathcal{{U}}_{2}}^{-1}$
is of the form
\begin{eqnarray}\label{equ531}
\left\{\begin{array}{ll}\tau_{1}=\tau+\alpha+\delta \Gamma(\varrho)+\delta\mathcal{V}_{1}\circ{\mathcal{{U}}_{2}}^{-1}(\tau,\varrho,\delta),\\
\varrho_{1}=\varrho+\delta\mathcal{V}_{2}\circ{\mathcal{{U}}_{2}}^{-1}(\tau,\varrho,\delta),
\end{array}\right.
\end{eqnarray}
where
$$\mathcal{V}_{1}=\left[\frac{\partial \tau}{\partial \theta}\widetilde{f}+\frac{\partial \tau}{\partial \rho}\widetilde{g}\right]+O_{1}(\delta)
,\ \ \ \ \ \ \mathcal{V}_{2}=\left[\frac{\partial I}{\partial \theta}\widetilde{f}+\frac{\partial I}{\partial \rho}\widetilde{g}\right]+O_{2}(\delta)
,$$
and the remainder term $O_{1}(\delta)$ is composed by $\delta, \widetilde{f}, \widetilde{g}$ and second-order derivatives of $\tau$,
the remainder term $O_{2}(\delta)$ is composed by second-order derivatives of $I$ and $\delta, \widetilde{f}, \widetilde{g}.$
Hence, the remainder terms $O_{1}(\delta), O_{2}(\delta)$ are real analytic and satisfy $O_{1}(\delta), O_{2}(\delta)\rightarrow0$ as $\delta\rightarrow0.$

Since $I(\theta,\rho)$ and $\widehat{L}(\theta,\varrho(\theta,\rho))$ are even in $\theta$, one has that the transformed mapping (\ref{equ531}) is reversible with respect to the involution $(\tau,\varrho)\mapsto(-\tau,\varrho).$

As we did in the proof of Theorem \ref{Th3}, we can obtain the existence of invariant curves of $\mathfrak{M}_{\delta}$ from Theorem \ref{Th2}.
This end the proof of Theorem \ref{Th4}.\qed
\section{Application}
In this section, we will apply the above results to the following nonlinear oscillator
\begin{eqnarray}\label{equ61}
x''+g(x)x'+\varpi^{2}x+\varphi(x)=f(t),
\end{eqnarray}
where $\varpi>$, $g,\varphi$ and $f$ are odd functions, $f(t)$ is a real analytic almost periodic function with the frequency $\omega=(\cdots,\omega_{\lambda},\cdots)$ and admits a spatial series expansion of the type (\ref{equ21}). We suppose that
\begin{eqnarray}\label{equ62}
\lim_{x\rightarrow+\infty}\varphi(x)=\varphi(+\infty)\in\mathbb{R},\ \ \ \ \ \ \ \lim_{|x|\rightarrow+\infty}x^{k}\varphi^{(k)}(x)=0, \ \ \ \ k\rightarrow\infty,
\end{eqnarray}
and
\begin{eqnarray}\label{equ63}
|x^{k}G^{(k)}(x)|\leq M , \ \ \ \ x\in\mathbb{R},\ \ \ k\geq0,
\end{eqnarray}
for some constant $M>0$, where $ G(x)=\int_{0}^{x}g(\xi)d\xi.$

The equation (\ref{equ61}) is equivalent to the system
\begin{eqnarray}\label{equ64}
\left\{\begin{array}{ll}x'=\varpi y-G(x),\\
y'=-\varpi x-\varpi^{-1}\varphi(x)+\varpi^{-1}f(t).
\end{array}\right.
\end{eqnarray}
Then the system (\ref{equ64}) is reversible with respect to the involution $(x,y)\mapsto(x,-y).$

Introduce the following  polar coordinates $(r, \theta)\mapsto(x, y)$:
$$x=r\sin\theta,\ \ \ \ \ \ \ \ y=r\cos\theta,$$
the system (\ref{equ64}) is changed into the form
\begin{eqnarray}\label{equ65}
\left\{\begin{array}{ll}r'=\varpi^{-1}[f(t)-\varphi(r\sin\theta)]\cos\theta-G(r\sin\theta)\sin\theta,\\
\theta'=\varpi-\varpi^{-1}r^{-1}[f(t)-\varphi(r\sin\theta)]\sin\theta-r^{-1}G(r\sin\theta)\cos\theta.
\end{array}\right.
\end{eqnarray}
According to our assumptions,
$$\theta'(t)\geq 0$$
for $r\gg1.$ This means that $t\mapsto \theta(t)$ is globally invertible. Denoting by $\theta\mapsto t(\theta)$ the inverse function, one has $\theta\mapsto (r(t(\theta), t(\theta)))$ solves
\begin{eqnarray}\label{equ66}
\left\{\begin{array}{lll}\displaystyle\frac{d r}{d\theta}=p_{1}(r,t,\theta),\\
\displaystyle\frac{d t}{d\theta}=p_{2}(r,t,\theta),
\end{array}\right.
\end{eqnarray}
where $r(\theta)=r(t(\theta),$ and
$$p_{1}(r,t,\theta)=\frac{\varpi^{-1}[f(t)-\varphi(r\sin\theta)]\cos\theta-G(r\sin\theta)\sin\theta}{\varpi-\varpi^{-1}r^{-1}[f(t)-\varphi(r\sin\theta)]\sin\theta-r^{-1}G(r\sin\theta)\cos\theta},$$
$$p_{2}(r,t,\theta)=\frac{1}{\varpi-\varpi^{-1}r^{-1}[f(t)-\varphi(r\sin\theta)]\sin\theta-r^{-1}G(r\sin\theta)\cos\theta}.$$
The system (\ref{equ66}) is reversible with respect to $(r,t)\mapsto(r,-t)$ since $p_{1}(r,-t,-\theta)=-p_{1}(r,t,\theta)$ and $p_{2}(r,-t,-\theta)=p_{2}(r,t,\theta)$ according to our symmetry assumptions.
Moreover the system (\ref{equ66}) can be written in the form
\begin{eqnarray}\label{equ67}
\left\{\begin{array}{lll}\displaystyle\frac{d r}{d\theta}=\varpi^{-2}[f(t)-\varphi(r\sin\theta)]\cos\theta-\varpi^{-1}G(r\sin\theta)\sin\theta+O(r^{-1}),\\
\displaystyle\frac{d t}{d\theta}=\varpi^{-1}-\varpi^{-3}r^{-1}[f(t)-\varphi(r\sin\theta)]\sin\theta-\varpi^{-2}r^{-1}G(r\sin\theta)\cos\theta+O(r^{-2}).
\end{array}\right.
\end{eqnarray}

In the following, we need to transform (\ref{equ67}) further such that we can use the invariant theorem to prove the existence of invariant curves for the Poincar\'{e} map of (\ref{equ67}).

Similar to \cite{KKL}, we have the following lemmas.
\begin{lemma}\label{lemma61}
Under the transformation $(r,t)\mapsto (\varrho,t),$ where
$$\varrho=r+S(r,\theta), \ \ \ \ S(r,\theta)=\int_{0}^{\theta}[\varpi^{-2}\varphi(r\sin\xi)\cos\xi+\varpi^{-1}G(r\sin\xi)\sin\xi]d\xi,$$
the system (\ref{equ67}) is transformed into
\begin{eqnarray}\label{equ68}
\left\{\begin{array}{lll}\displaystyle\frac{d \varrho}{d\theta}=\varpi^{-2}f(t)\cos\theta+O(\varrho^{-1}),\\
\displaystyle\frac{d t}{d\theta}=\varpi^{-1}-\varpi^{-3}\varrho^{-1}[f(t)-\varphi(\varrho\sin\theta)]\sin\theta-\varpi^{-2}\varrho^{-1}G(\varrho\sin\theta)\cos\theta+O(\varrho^{-2}).
\end{array}\right.
\end{eqnarray}
\end{lemma}

The system (\ref{equ68}) is reversible with respect to $(\varrho,t)\mapsto(\varrho,-t)$ by the properties of reversible systems.
\begin{lemma}
For $k\geq0$ and $J(\varrho)=\frac{1}{2\pi\varrho}\int_{0}^{2\pi}\varphi(\varrho\sin\xi)\sin\xi d\xi,$ one has
\begin{eqnarray}\label{equ69}
\lim_{\varrho\rightarrow+\infty}\varrho^{k+1}J^{(k)}(\varrho)=(-1)^{k}k!\frac{2}{\pi}\varphi(+\infty).
\end{eqnarray}
\end{lemma}

\begin{lemma}
Under the transformation $(\varrho,t)\mapsto (\varrho,\tau),$ where
$$\tau=t+T(\varrho,\theta), \ \ \ \ T(\varrho,\theta)=-\int_{0}^{\theta}\varpi^{-3}\varrho^{-1}[\varphi(\varrho\sin\xi)\sin\xi-\varrho J(\varrho)]d\xi+\varpi^{-2}\varrho^{-1}\int_{0}^{\theta}G(\varrho\sin\xi)\cos\xi d\xi,$$
the system (\ref{equ68}) is transformed into
\begin{eqnarray}\label{equ610}
\left\{\begin{array}{lll}\displaystyle\frac{d \varrho}{d\theta}=\varpi^{-2}f(\tau)\cos\theta+O(\varrho^{-1}),\\
\displaystyle\frac{d \tau}{d\theta}=\varpi^{-1}+\varpi^{-3}J(\varrho)-\varpi^{-3}\varrho^{-1}f(\tau)\sin\theta+O(\varrho^{-2}).
\end{array}\right.
\end{eqnarray}
\end{lemma}

\begin{lemma}
For \begin{eqnarray}\label{equ611}
\eta(\varrho)=J(\varrho)-\frac{2}{\pi}\varphi(+\infty)\varrho^{-1},
\end{eqnarray}
\end{lemma}
one has
\begin{eqnarray}\label{equ612}
\eta(\varrho)=o(\varrho^{-1}).
\end{eqnarray}

By (\ref{equ611}) and (\ref{equ612}),  one has
\begin{eqnarray}\label{equ613}
\left\{\begin{array}{lll}\displaystyle\frac{d \varrho}{d\theta}=\varpi^{-2}f(\tau)\cos\theta+O(\varrho^{-1}),\\
\displaystyle\frac{d \tau}{d\theta}=\varpi^{-1}+\frac{2}{\pi}\varpi^{-3}\varphi(+\infty)\varrho^{-1}-\varpi^{-3}\varrho^{-1}f(\tau)\sin\theta+o(\varrho^{-1}).
\end{array}\right.
\end{eqnarray}
We can verify that this system is reversible with respect to $(\varrho,\tau)\mapsto(\varrho,-\tau).$

Introducing a new variable $\rho$ and a small parameter $\epsilon>0$ by $$\varrho^{-1}=\epsilon\rho,\ \ \  \rho\in[1,2].$$
Obviously, $\varrho\gg1\Leftrightarrow\epsilon\ll1.$ The system (\ref{equ613}) is changed to the form
\begin{eqnarray}\label{equ614}
\left\{\begin{array}{lll}\displaystyle\frac{d \rho}{d\theta}=-\epsilon\rho^{2}\varpi^{-2}f(\tau)\cos\theta+O(\epsilon^{2}),\\
\displaystyle\frac{d \tau}{d\theta}=\varpi^{-1}+\frac{2}{\pi}\epsilon\rho\varpi^{-3}\varphi(+\infty)-\varpi^{-3}\epsilon\rho f(\tau)\sin\theta+O(\epsilon^{2}).
\end{array}\right.
\end{eqnarray}
The system (\ref{equ614}) is reversible with respect to $(\rho,\tau)\mapsto(\rho,-\tau)$ since the transformation $(\varrho,\tau)\mapsto(\rho,\tau)$ is reversible with respect to $(\varrho,\tau)\mapsto(\varrho,-\tau).$

Suppose that the solution $(\rho(\theta),\tau(\theta))$ of (\ref{equ614}) has the following expression
\begin{eqnarray}\label{equ615}
\rho(\theta;\rho_{0},\tau_{0})=\rho_{0}+\epsilon F_{1}(\theta;\rho_{0},\tau_{0},\epsilon),\ \ \ \tau(\theta;\rho_{0},\tau_{0})=\tau_{0}+\varpi^{-1}\theta+\epsilon F_{2}(\theta;\rho_{0},\tau_{0},\epsilon),
\end{eqnarray}
that is, one has
\begin{eqnarray}\label{equ616}
F_{1}(0;\rho_{0},\tau_{0},\epsilon)=F_{2}(0;\rho_{0},\tau_{0},\epsilon)=0.
\end{eqnarray}

Differentiating (\ref{equ615}) and comparing with (\ref{equ614}), one has
\begin{eqnarray*}
\left\{\begin{array}{llll}\displaystyle\frac{\partial F_{1}}{\partial\theta}=-\varpi^{-2}(\rho_{0}+\epsilon F_{1})^{2}f(\tau_{0}+\varpi^{-1}\theta+\epsilon F_{2})\cos\theta+ O(\epsilon),\\
\displaystyle\frac{\partial F_{2}}{\partial\theta}=\frac{2}{\pi}\varpi^{-3}(\rho_{0}+\epsilon F_{1})\varphi(+\infty)-\varpi^{-3}(\rho_{0}+\epsilon F_{1}) f(\tau_{0}+\varpi^{-1}\theta+\epsilon F_{2})\sin\theta+ O(\epsilon).
\end{array}\right.
\end{eqnarray*}
Hence
\begin{eqnarray}\label{equ617}
\left\{\begin{array}{llll}F_{1}(2\pi;\rho_{0},\tau_{0},\epsilon)=-\varpi^{-2}\rho_{0}^{2}\int_{0}^{2\pi}f(\tau_{0}+\varpi^{-1}\theta)\cos\theta d\theta+ O(\epsilon),\\
F_{2}(2\pi;\rho_{0},\tau_{0},\epsilon)=\varpi^{-3}\rho_{0}\left(4\varphi(+\infty)-\int_{0}^{2\pi}f(\tau_{0}+\varpi^{-1}\theta)\sin\theta d\theta\right)+ O(\epsilon).
\end{array}\right.
\end{eqnarray}

The Poincar\'{e} map $P_{2\pi}$ of the system (\ref{equ614})
has the expansion
\begin{eqnarray}\label{equ618}
P_{2\pi}:\left\{\begin{array}{ll}\rho_{1}=\rho_{0}+\epsilon m(\rho_{0},\tau_{0})+O(\epsilon^{2})
,\\ \tau_{1}=\tau_{0}+\varpi^{-1}2\pi+\epsilon l(\rho_{0},\tau_{0})+O(\epsilon^{2}),
\end{array}\right.
\end{eqnarray}
where
\begin{eqnarray*}
m(\rho_{0},\tau_{0})&=& -\rho_{0}^{2}\varpi^{-2}\int_{0}^{2\pi}f(\tau_{0}+\varpi^{-1}\theta)\cos\theta d\theta= \rho_{0}^{2}\varpi^{-3}\int_{0}^{2\pi}f'(\tau_{0}+\varpi^{-1}\theta)\sin\theta d\theta,
\\l(\rho_{0},\tau_{0})&=& \rho_{0}\varpi^{-3}\left(4\varphi(+\infty)-\int_{0}^{2\pi}f(\tau_{0}+\varpi^{-1}\theta)\sin\theta d\theta\right).
\end{eqnarray*}
We see that $P_{2\pi}$ is reversible with respect to the involution $(\rho_{0},\tau_{0})\mapsto(\rho_{0},-\tau_{0}).$

By Theorem \ref{Th3} and Theorem \ref{Th4}, one may obtain the following conclusion.

\begin{theorem}\label{Th61} In the previous settings, we assume that
$\omega_{A}=\{\omega_{\lambda}:\lambda\in A\}$ and $\varpi$ are rationally independent for any $A\in\mathcal{S}$, and
$$\varphi(+\infty)\neq0.$$
Then system (\ref{equ61}) has many almost periodic solutions and all the solutions are bounded.
\end{theorem}　
{\bf Proof}. From Theorem \ref{Th3}, we know that $P_{2\pi}$
has invariant curves in the domain $\mathbb{R}\times[1,2]$ if $\epsilon$ is sufficiently small and
\begin{eqnarray*}\lim_{T\rightarrow\infty}\frac{1}{T}\int_{0}^{T}l(\rho_{0},\tau_{0})d\tau_{0}\neq0.
\end{eqnarray*}
By (\ref{equ21}), $f$ has the following series expansion
\begin{eqnarray}\label{equ}
f(t)=f_{0}+\sum_{A\in\mathcal{S}}\sum_{\substack{k\neq0\\ {\rm supp}k\subseteq A}}f_{k}e^{i\langle k,\omega\rangle t}.
\end{eqnarray}
By Fubini's theorem, it follows that
\begin{eqnarray}
\lim_{T\rightarrow\infty}\frac{1}{T}\int_{0}^{T}l(\rho_{0},\tau_{0})d\tau_{0}=-\varpi^{-3}4\varphi(+\infty)
\end{eqnarray}
Hence if $\varphi(+\infty)\neq0$, then the existence of invariant curves as well as the boundedness of solutions are guaranteed by Theorem 5.4.\qed

\begin{theorem}\label{Th62}
If there is some set $A\in\mathcal{S}$ such that $\omega_{A}=\{\omega_{\lambda}:\lambda\in A\}$ and $\varpi$ are rationally dependent. Denote by $\mathfrak{S}$ the set of all $A\in \mathcal{S} $ such that there is the integer vector
$k\in\mathbb{Z}^{\mathbb{Z}}_{\mathcal{S}}\setminus\{0\}$ and $suppk\subseteq A, \langle k,\omega\rangle\varpi^{-1}\in\mathbb{Z}.$
We denote by $f_{\mathfrak{S}}$ the function
\begin{eqnarray}
f_{\mathfrak{S}}=\sum_{A\in\mathfrak{S}}\sum_{{\rm supp}k\subseteq A}f_{k}e^{i\langle k,\omega\rangle t}.
\end{eqnarray}
If
\begin{eqnarray*}
-4\varphi(+\infty)+\int_{0}^{2\pi}f_{\mathfrak{S}}(\tau_{0}+\varpi^{-1}\theta)\sin\theta d\theta\neq0,\ \ \ \ \forall \tau_{0}\in\mathbb{R}.
\end{eqnarray*}
Then system (\ref{equ61}) has many almost periodic solutions and all the solutions are bounded.
\end{theorem}　
{\bf Proof}. Without loss of generality, we assume
\begin{eqnarray*}
-4\varphi(+\infty)+\int_{0}^{2\pi}f_{\mathfrak{S}}(\tau_{0}+\varpi^{-1}\theta)\sin\theta d\theta>0,\ \ \ \ \forall \tau_{0}\in\mathbb{R}.
\end{eqnarray*}
Then we choose the function $I$ as
\begin{eqnarray*}
I(\rho_{0},\tau_{0})=\frac{\rho_{0}}{-4\varphi(+\infty)+\int_{0}^{2\pi}f_{\mathfrak{S}}(t+\varpi^{-1}\theta)\sin\theta d\theta}.
\end{eqnarray*}
We can verify that the assumptions in Theorem \ref{Th4} are satisfied. Hence, the existence of invariant curves and the boundedness of solutions are proved.\qed

\noindent{$\mathbf{Acknowledgments}$}

{\indent This work was supported by the NSFC (grant no. 11571327).}

\vskip1cm


\end{spacing}

\end{document}